\documentclass[11pt]{amsart}

\usepackage{amssymb,upgreek}

\usepackage{url}
\usepackage{bm}
\usepackage[left=1in,top=1in,right=1in,bottom=1in,head=.2in]{geometry}
\usepackage[titletoc]{appendix}

\usepackage[textsize=scriptsize]{todonotes}
\usepackage{color}
\usepackage{mathtools} 

\usepackage{fancyhdr}
\usepackage{mathrsfs}
\usepackage[cmtip,all]{xy}

\usepackage{hyperref}
\usepackage{enumerate}

\swapnumbers
\newtheorem{theorem}{Theorem}[section] 
\newtheorem{lemma}[theorem]{Lemma}

\newtheorem*{theorem*}{Theorem}

\newtheorem*{fcthm*}{Finite Cork Theorem}
\newtheorem*{ccthm*}{Cork Consolidation Theorem}
\newtheorem*{thm1*}{Theorem 1}
\newtheorem*{thm2*}{Theorem 2}
\newtheorem*{lbthm*}{Generalized 4D Lightbulb Theorem}
\newtheorem*{2lbthm*}{Generalized 4D Lightbulb Theorem (restated)}
\newtheorem*{icthms*}{Infinite Cork Theorems}
\newtheorem*{aclemma*}{\ac-Lemma}
\newtheorem*{mclemma*}{Multicork Lemma}
\newtheorem*{multicorktheorem*}{Multicork Theorem}
\newtheorem*{lemma*}{Lemma}
\newtheorem*{corollary*}{Corollary}

\newcommand{\thistheoremname}{}
\newtheorem{genericthm}[theorem]{\thistheoremname}

\theoremstyle{definition}

\newtheorem{remark}[theorem]{Remark}

\newtheorem*{remark*}{Remark}
\newtheorem*{definition*}{Definition}
\newtheorem*{remarks*}{Remarks}
\newtheorem*{addenda*}{Addenda}

\newcommand{\pf}{\vskip-5pt \vskip-5pt \proof}

\newcommand{\fig}[3]{\begin{figure}\includegraphics[height=#1pt]{#2}#3\end{figure}}

\newcommand{\bit}[1]{\textbf{\textit{#1}}} 

\newcommand{\bc}{\mathbb C}

\newcommand{\id}{\textup{id}}
\newcommand{\pt}{\textup{pt}}

\newcommand{\sto}{\!\!\xymatrix@C=1em{{}\ar@{~>}[r]&{}}\!\!}

\newcommand{\sss}{S^2\!\!\times\!S^2}

\newcommand{\cptwo}{\bc P^2}

\newcommand{\interior}{\textup{int}}
\newcommand{\ac}{\textup{AC}}

\newcommand{\lk}{\textup{lk}}

\newcommand{\cs}{\mathop\#}

\newcommand{\foot}[1]{\setcounter{footnote}{1}\footnote{\ #1}}

\newcommand{\items}{\begin{itemize}[leftmargin=25pt,rightmargin=5pt]
  \setlength\itemsep{2pt}}
\newcommand{\stopitems}{\end{itemize}}

\address{Bryn Mawr College,
Bryn Mawr, PA 19003}
\email{hrschwartz@brynmawr.edu} 

\begin{document}

\title{Equivalent non-isotopic spheres in $4$-manifolds}
\author{Hannah R. Schwartz}

\begin{abstract}
We construct infinitely many smooth oriented  $4$-manifolds containing pairs of homotopic, smoothly embedded $2$-spheres that are \emph{not} topologically isotopic, but that are equivalent by an ambient diffeomorphism inducing the identity on homology. These examples show that Gabai's recent ``Generalized" 4D Lightbulb Theorem does not generalize to arbitrary $4$-manifolds. In contrast, we also show that there are smoothly embedded $2$-spheres that are both equivalent and topologically isotopic, but not smoothly isotopic.
\end{abstract}

\maketitle

\vskip-.4in
\vskip-.4in

\parskip 2pt

\setcounter{section}{-1}

\parskip 2pt

\section{Introduction and Motivation}

The 3D Lightbulb Theorem gives an isotopy in $S^1 \times S^2$ from the curve $S^1 \times \{\pt\}$ to any embedded curve intersecting $\{\pt\} \times S^2$ transversally in exactly one point. A recent result due to Gabai \cite{dave:lightbulb} extends this theorem to $4$ dimensions, showing that any smoothly embedded $2$-sphere in $S^2 \times S^2$ homologous to $S^2 \times \{\pt\}$ and intersecting $\{\pt\} \times S^2$ transversally in one point  is smoothly isotopic to $S^2 \times \{\pt\}$.\foot{It has been known since the $80$'s that there exists a diffeomorphism pseudo-isotopic to the identity carrying $S$ to $S^2 \times \{\pt\}$ by \protect\cite{litherland}, and that the spheres $S$ and $S^2 \times \{\pt\}$ are topologically isotopic by \protect\cite{marumoto}.} The key feature in both theorems is the existence of an embedded $2$-sphere of square zero transversally intersecting the submanifold being isotoped in a single point, which we call a \bit{dual} to the submanifold. In fact, Gabai proves a similar result for a broader class of $4$-manifolds.

\begin{lbthm*}\label{thm1} 
Let $S$ and $T$ be homotopic $2$-spheres smoothly embedded in a smooth, orientable $4$-manifold $X$, with a common dual. If $\pi_1(X)$ has no $2$-torsion, then there is a smooth ambient isotopy from $S$ to $T$ in $X$.
\end{lbthm*} 

Although homotopy implies (locally-flat) \emph{topological} isotopy for surfaces in sufficiently nice settings (details can be found in \cite{sunukjian:isotopy}), this is not always the case for \emph{smooth} isotopy without the existence of a common dual. For instance, work of Donaldson \cite{donaldson} and Wall \cite{wall:4-manifolds} can be used to construct examples of smoothly embedded, homotopic 2-spheres with simply-connected complements that are not smoothly isotopic (see also \cite{akmr:isotopy}). These arguments obstruct smooth isotopy by producing pairs of non-diffeomorphic $4$-manifolds by either surgering or blowing down smoothly embedded pairs of homotopic spheres, and hence preclude even a diffeomorphism of the ambient manifold taking one sphere to the other.  

Even with a common dual, we show that there are instances in which homotopic embedded spheres fail to be topologically isotopic. However in this case, by virtue of their common dual, by Lemma \ref{diff} the spheres are \bit{equivalent}, meaning that there is an orientation preserving self-diffeomorphism of $X$ (restricting to the identity on $\partial X$ when it is non-empty) carrying one sphere to the other and inducing the identity map on $H_2(X; \mathbb{Z})$. 

\begin{thm1*}\label{thm2} 
There are infinitely many smooth $4$-manifolds $X$ with $\pi_1(X) \cong \mathbb{Z}_2$ containing smoothly embedded, homotopic $2$-sphere pairs that are equivalent but not topologically concordant in $X \times I$. 
\end{thm1*}

In particular, without the 2-torsion assumption or some other hypotheses, the Generalized 4D Lightbulb Theorem is false, giving a negative answer to Question $10.14$ posed in an early draft of \cite{dave:lightbulb}. The approach here is topologically rather than smoothly based, yet (unlike the constructions referenced above) allows the spheres to be related by a diffeomorphism of the ambient manifold. We also obtain the following, relying on results of Morgan and Szab\'o \cite{morganszabo} coming from gauge theory. 

\begin{thm2*}There exist simply-connected $4$-manifolds containing pairs of smoothly embedded $2$-spheres that are both equivalent and topologically isotopic, but not smoothly isotopic.
\end{thm2*}

Our results illustrate the curious distinction between \emph{smooth equivalence, topological isotopy, and smooth isotopy} of surfaces  in dimension $4$. In contrast to both theorems above, examples of spheres that are topologically isotopic but not smoothly equivalent (and hence not smoothly isotopic) can also be cooked up, by applying the aforementioned work of Donaldson \cite{donaldson} and Wall \cite{wall:4-manifolds}. 

\smallskip
\noindent
{\bf Acknowledgments.} I extend my warmest thanks to Dave Gabai for his guidance and kind encouragement during this project. Thank you also to my wonderful adviser Paul Melvin and academic ``baby brother" Isaac Craig for their never-ending patience and enlightening insight, to both Danny Ruberman and Dave Auckly for their helpful advice and suggestions, and to the referee for their careful edits. 

\section{Background}

\fig{300}{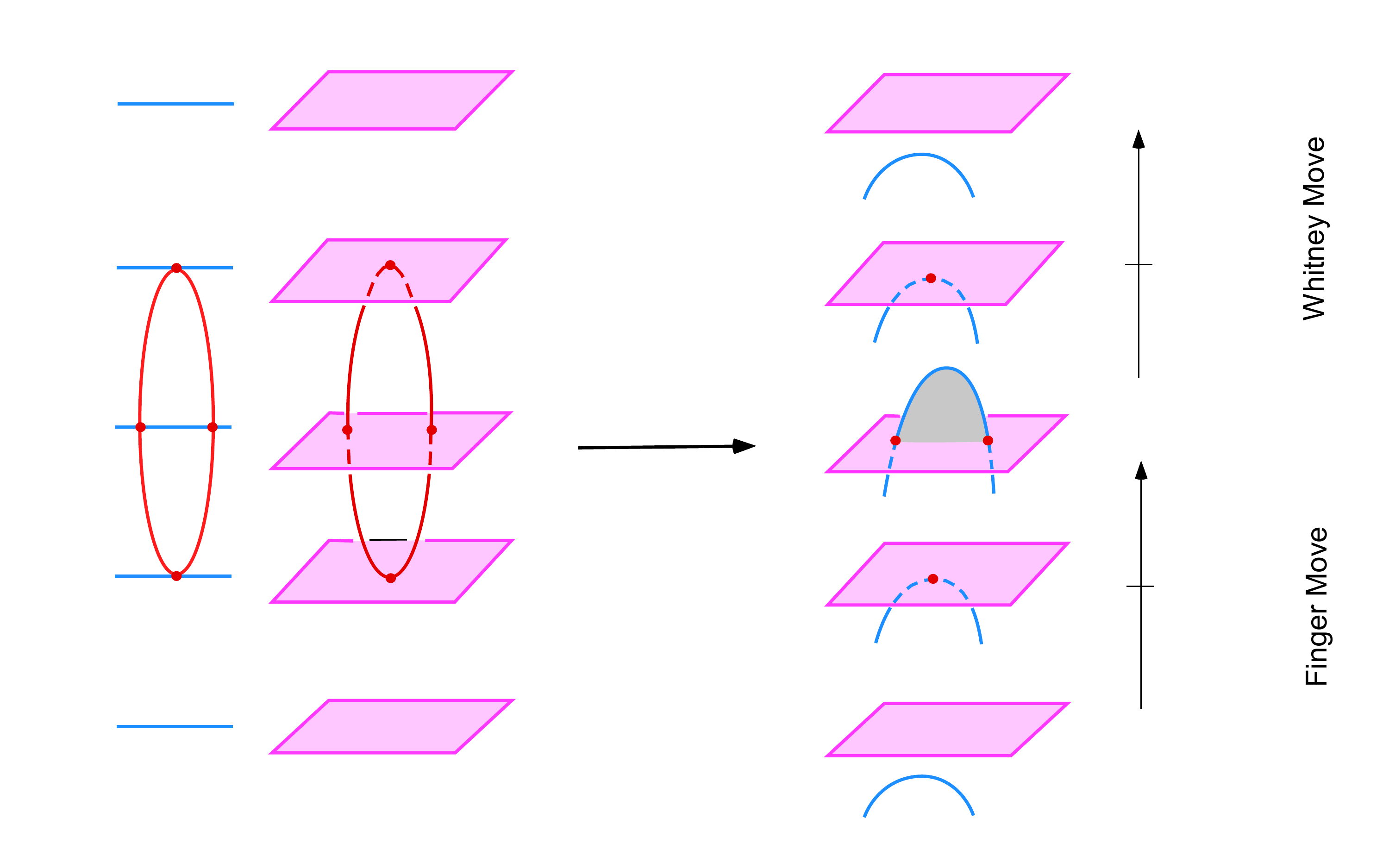}{
\put(-159,156){$W$}
\put(-181,142){\small $-$}
\put(-136,142){\small $+$}
\put(-255,150){$H$}
\put(-75,93){$t_{min}$}
\put(-75,204){$t_{max}$}
\caption{A finger move is performed along an embedded arc, whereas a Whitney move is supported in the regular neighborhood of an embedded ``Whitney disk" $W$. Given a map $H: S^2 \times I \to X \times I$ induced by a regular homotopy $h$, the cycles of $h$ have two minima with respect to the $I$ factor (drawn vertically) for each finger move, and two maxima for each Whitney move.}
\label{fw}}

By Smale \cite[Theorem D]{smale}, smoothly embedded $2$-spheres in a smooth, orientable $4$-manifold $X$ are homotopic if and only if they are \bit{regularly homotopic}, i.e. homotopic through smooth immersions. Regular homotopies of surfaces in $4$-dimensions are particularly convenient to analyze, since generically there are only finitely many times during the homotopy at which the immersed surface is not self-transverse -- this occurs when double points of opposite sign are either introduced or cancelled. The local model for the regular homotopy removing pairs of double points is called a \bit{Whitney move}, and its inverse is called a \bit{finger move}. Both of these homotopies are depicted in Figure \ref{fw} and Figure \ref{reghom1}.\foot{See \protect\cite{freedman-quinn:4-manifolds} and also Casson's lectures in \protect\cite{casson} for more exposition.}

\fig{350}{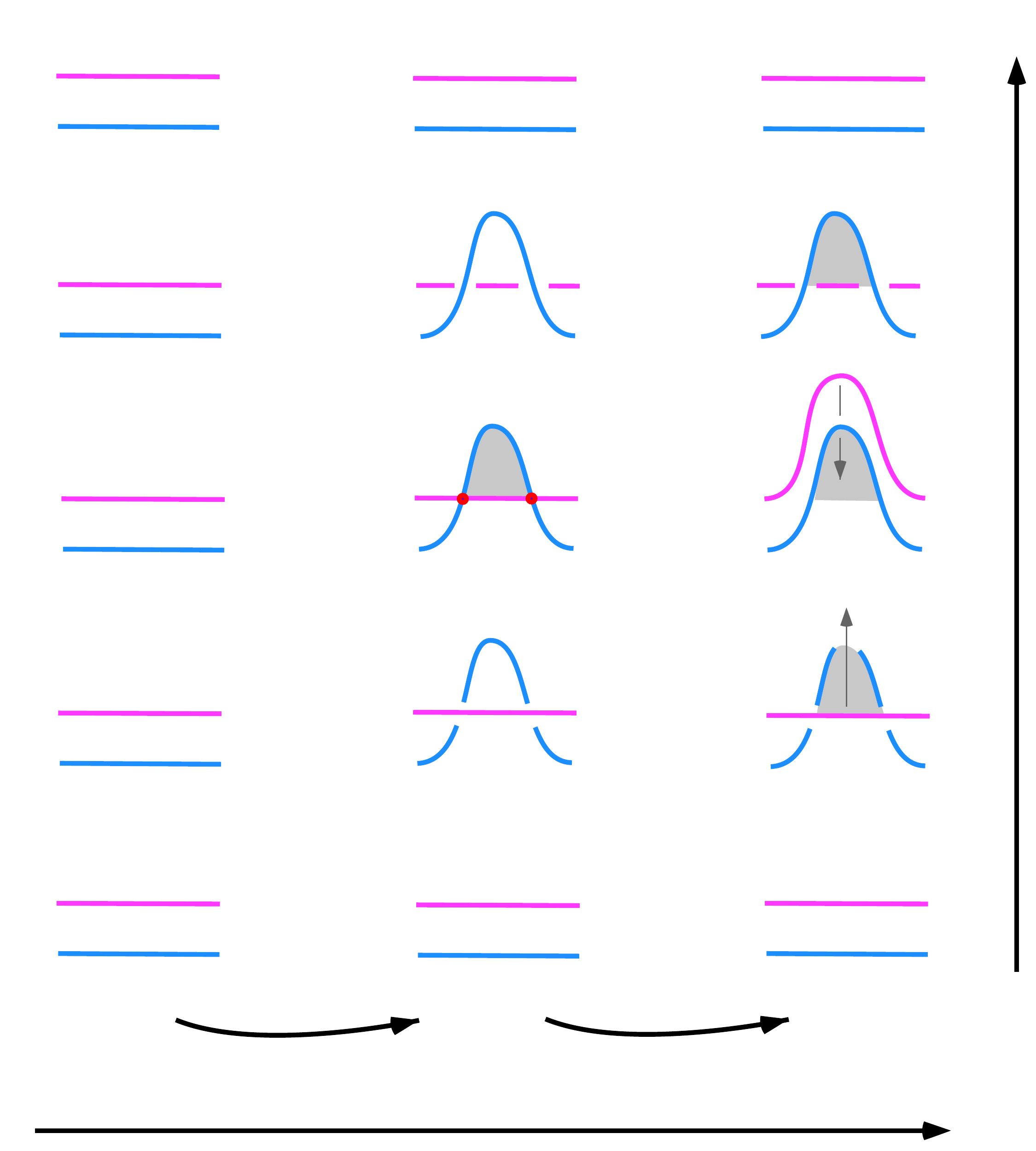}{
\put(-184,203){\small{$-$}}
\put(-152,203){\small{+}}
\put(-169,202){\small{$W$}}
\put(-254,20){\small{Finger move}}
\put(-145,20){\small{Whitney move}}
\put(-240,-10){\small{Time in the regular homotopy}}
\caption{The local model in $\mathbb{R}^3 \times \mathbb{R}$ of an immersed sphere both before and after a finger and Whitney move. Each level set, or ``slice", is the intersection of the sphere with $\mathbb{R}^3 \times \{h\}$ for some height $h \in \mathbb{R}$ recorded by the vertical axis. Note that after the Whitney move using the Whitney disk $W$, the cancelling double points are avoided by pushing one sheet of the sphere first along the ``front" and then the ``back" of the Whitney disk $W$.}
\label{reghom1}}

Suppose $h: S^2 \times I \to X$ is a regular homotopy between $2$-spheres $S$ and $T$. Then the map $H:S^2 \times I \to X \times I$ sending $(x,t) \mapsto (h(x,t),t)$ is a smooth immersion with preimage of its double points a union of disjointly embedded curves in $S^2 \times I$, which we call the \bit{cycles} of the regular homotopy. The cycles can be isotoped so that they are transverse to the level sets $S^2 \times \{t\}$ except at the times $t_{min}, t_{max} \in I$ where they attain maxima and minima with respect to the $I$ factor. Each minimum (resp. maximum) is a point in the preimage of some double point at which the immersion is non-transverse during a finger (resp. Whitney) move. 

A cycle which non-trivially double covers its image is called a \bit{crossed cycle}. An example is shown in Figure \ref{fw2}. A crossed cycle is said to \bit{correspond} to the element of $\pi_1(X)$ represented by the loop in $X \times I$ that it double covers, projected to $X$. This element is well-defined (choosing the basepoint of $\pi_1(X)$ in the image of $h$) since free homotopy is equivalent to homotopy in $S^2 \times I$.  

\fig{200}{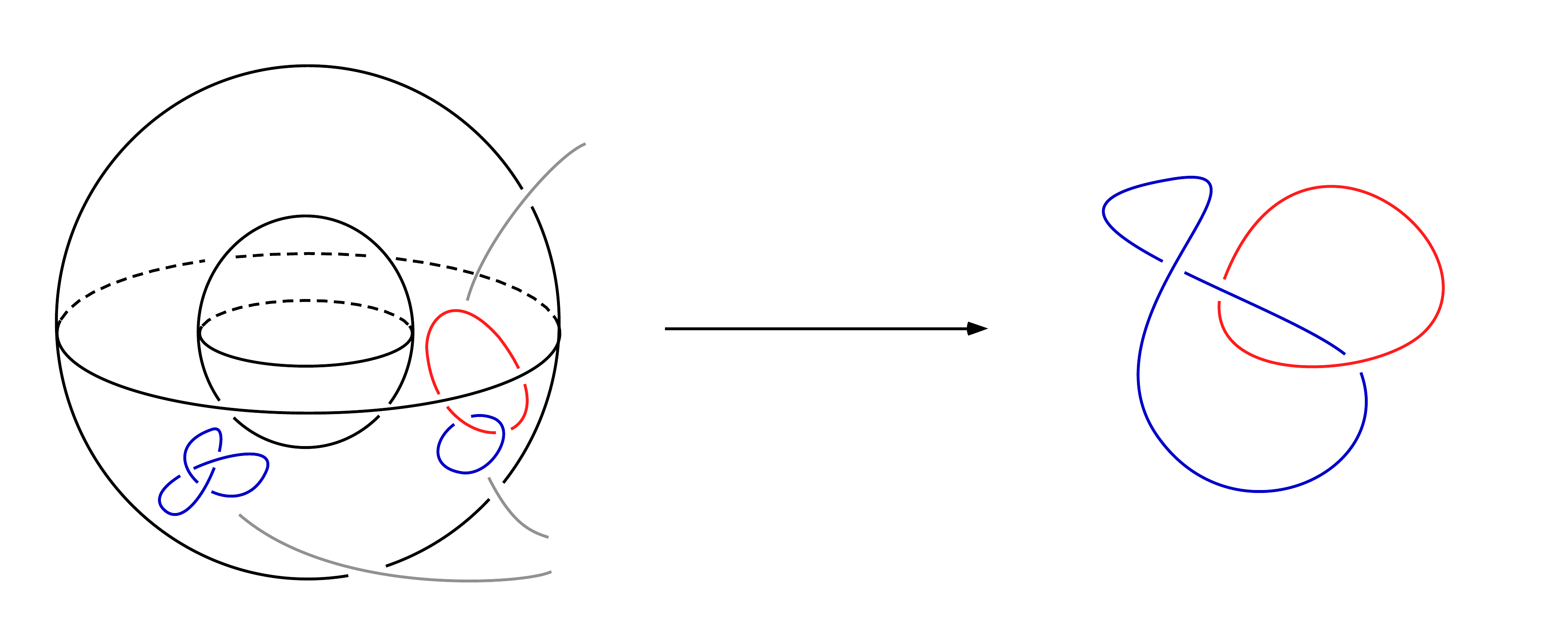}{
\put(-235,105){$H$}
\put(-300,154){Crossed cycle}
\put(-310,24){Uncrossed cycles}
\caption{The double points (shown on the right) of a map $H: S^2 \times I \to X \times I$ induced by a regular homotopy, and their preimages (shown on the left).}
\label{fw2}}

Note that the crossed cycles of a homotopy necessarily correspond to elements of order at most 2, and give rise to \emph{double tubes} (discussed in Sections 5 and 6 of \cite{dave:lightbulb}) that represent these elements. Hence, Gabai's Theorem 1.3 in \cite{dave:lightbulb} translates to the following.

\begin{2lbthm*} 
Let $S$ and $T$ be smoothly embedded 2-spheres in a smooth, orientable $4$-manifold $X$. Suppose that $S$ and $T$ 
\begin{enumerate}
\item [\small\bf 1)] have a common dual, and 
\item [\small\bf 2)] are related by a regular homotopy supported away from a neighborhood of the dual, with an even number of crossed cycles corresponding to each element of order 2 in $\pi_1(X)$. 
\end{enumerate}
Then, $S$ and $T$ are smoothly isotopic in $X$. 
\end{2lbthm*}

We will show that condition $\small\bf 2)$ is necessary, by constructing examples of homotopic $2$-spheres with a common dual that are \emph{not} smoothly (or even topologically) isotopic. 
\section{Main Examples} 

Consider the family of smooth, orientable, compact 4-manifolds $X_{p,q}$ illustrated in Figure \ref{mfld}.\foot{We use the dotted circle notation for 1-handles as in \protect\cite[Chapter I.2]{kirby:4-manifolds}.} For each $p, q \in \mathbb{Z}$, $H_1(X_{p,q}) \cong \pi_1(X_{p,q}) \cong \mathbb{Z}_2$. The diffeomorphism type of $X_{p,q}$ is determined by the integer $p$, and the value of $q \mod 2$. To see this, first note that the manifolds $X_{p_1,q_1}$ and $X_{p_2,q_2}$ are diffeomorphic when $p_1 = p_2$ and $q_1 \equiv q_2 \mod 2$, as sliding the 2-handle with framing $q$ over the 2-handle with framing $0$ sufficiently many times reduces its framing to either $0$ or $1$. The manifolds are distinguished otherwise by their universal (2-fold) covers $\widetilde X_{p,q}$ since $H_1(\partial \widetilde X_{p,q})$ is determined by the value of $p$ (an exercise in Kirby calculus). 

\fig{130}{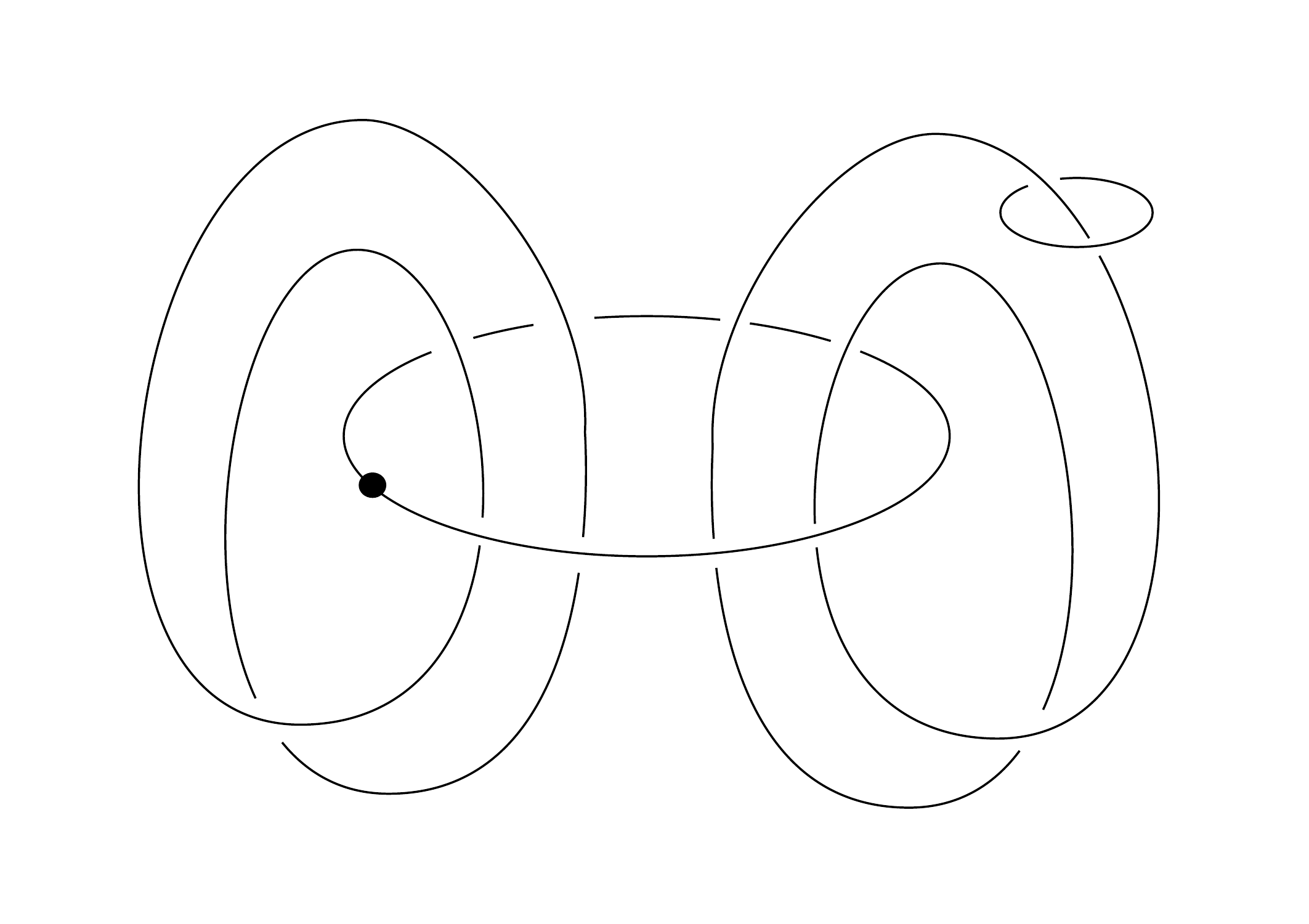}{
\put(-175,70){$p$}
\put(-16,66){$q$}
\put(-18,97){$0$}
\caption{The manifold $X_{p,q}$}
\label{mfld}}

From now on, let $X= X_{p,q}$ for any choice of framings $p,q \in \mathbb{Z}$. The union of the $0$-handle and the $1$-handle of $X$ is a copy of $S^1 \times D^2 \times [-1,1]$ which we label $X^{(1)}$. Letting $\gamma$ denote the (2,1) torus curve pushed into the interior of the solid torus $S^1 \times D^2$, the $2$-handles of framings $p$ and $q$ are attached along neighborhoods of the circles $\gamma_\pm=\gamma \times \{\pm 1\}$. 

Let $\Sigma_0$ denote the embedded 2-sphere in $X$ consisting of the cores of these $2$-handles connected by the annulus $\gamma \times I \subset X^{(1)}$. For each $n \in \mathbb{Z}$, let $\Sigma_n$ denote the image of $\Sigma_0$ under the automorphism $\phi_n: X \to X$ twisting the $S^1 \times D^2$ factor of $X^{(1)}$ longitudinally (i.e. in the $S^1$ direction) $n$ times as it descends through the $I$ factor, and fixing the $2$-handles of $X$. All $\Sigma_n$ are homologous of square $p+q$, and share a dual (the union of the core of the zero framed $2$-handle and the spanning disk in $S^3$ of its attaching curve). 

\begin{theorem}\label{spherediff}
If $i \equiv j \mod 4$, then the $2$-spheres $\Sigma_i$ and $\Sigma_j$ are smoothly isotopic in $X$.   
\end{theorem}

\pf  This follows as a corollary of the 4D Lightbulb Theorem. For, the lifts of these 2-spheres to the universal $2$-fold cover $\widetilde X$ are homologous, and hence homotopic, if and only if $i$ and $j$ have the same parity. So, there are regular homotopies $h:S^2 \times I \to X \times I$ from $\Sigma_i$ to $\Sigma_{i+2}$, and $h'=(\phi_2 \times \id) \circ h$ from $\Sigma_{i+2}$ to $\Sigma_{i+4}$, with the same number of crossed cycles corresponding to the $2$-torsion element generating $\pi_1(X) \cong \mathbb{Z}_2$. Performing first $h$ and then $h'$ therefore gives a regular homotopy from  $\Sigma_{i}$ to $\Sigma_{i+4}$ with an even number of crossed cycles corresponding to the generator of $\pi_1(X)$. Thus, by the restated version of the 4D Lightbulb Theorem, the 2-spheres $\Sigma_i$ and $\Sigma_{i+4}$ are smoothly ambiently isotopic in $X$.  \qed 

\begin{remark}\label{iso1}
In Appendix \ref{appendix}, we present an explicit regular homotopy from $\Sigma_i$ to $\Sigma_{i+2}$ with a single crossed cycle corresponding to the $2$-torsion element of $\pi_1(X)$, consisting of a single finger move and Whitney move. In some sense, our example is the smallest manifold in which such a homotopy is possible: the manifold $X$ includes only the support of this regular homotopy, and a neighborhood of the dual $2$-sphere. 
\end{remark}

The presence of a dual sphere also guarantees that there is a self-diffeomorphism of $X$ carrying $\Sigma_0$ to $\Sigma_n$ for all $n$, that (unlike $\phi_n$) restricts to the identity on $\partial X$. 

\begin{lemma} \label{diff}
Let $S$ and $S'$ be 2-spheres of the same square smoothly embedded in a smooth, orientable 4-manifold $X$, with a common dual $S^*$. Then, there is an orientation preserving diffeomorphism of pairs $\phi: (X, S) \to (X, S')$. Furthermore, if $\partial X$ is non-empty, then $\phi$ can be chosen to restrict to the identity on $\partial X$. 
\end{lemma}  

\emph{Proof.} Since $S$ and $S'$ have equal square, the regular neighborhoods $N,N' \subset X$ of $S \cup S^*$ and $S' \cup S^*$ can be identified by a diffeomorphism $h: N \to N'$ mapping $S$ to $S'$. Let $X^*$ denote the manifold obtained from $X$ by surgering $S^*$. This surgery restricted to $N$ and $N'$ results in 4-balls $B, B' \subset X^*$. By Palais' theorem \cite{palais}, there is a smooth ambient isotopy $F_t: X^* \to X^*$ carrying $B$ to $B'$. 

The end of this isotopy gives a diffeomorphism $F_1 : X^*-B \to X^*-B'$. By Cerf \cite{cerf}, the map $F_1$ restricted to $\partial B=S^3$ may be isotoped to $h |_{\partial N=\partial B}$ so that $F_1$ can be extended by across $N$ by $h$ to a diffeomorphism $\widehat F_1 : X \to X$ sending $S$ to $S'$. If $\widehat F_1$ does not restrict to the identity on $\partial X$, modify the diffeomorphism on a collar of $\partial X$ by inserting the isotopy $F|_{\partial X}$ from $F_1|_{\partial X}$ to $\id_{\partial X}$. This produces a smooth map $(X, S) \to (X, S')$ restricting to the identity on $\partial X$, as desired. 
\qed
\smallbreak

It follows from Lemma \ref{diff} that most standard methods of distinguishing smooth isotopy classes of 2-spheres fail when there is a common dual. For instance, when $S$ and $S'$ have trivial normal bundles, they \emph{cannot} be distinguished by the manifolds that result from $X$ by Gluck twisting \cite{gluck:2-spheres} or surgering either 2-sphere. Likewise, when $S$ and $S'$ have self-intersection $\pm 1$, they cannot be distinguished by blowing down. Note that since $F_t(B)$ is not necessarily equal to $B'$ for all $t \in I$, it is not possible to surger $X^*$ back to $X$ and then extend the isotopy along the resulting $D^2 \times S^2$  (as we did not assume the 2-spheres are homotopic or even homologous, we should not expect this). 

\begin{remark}
A similar result holds for homotopic $\Sigma^*$-inessential\footnote{The term ``$\Sigma^*$-inessential" is defined in \cite[Theorem $9.7$]{dave:lightbulb}.} surfaces $\Sigma, \Sigma' \subset X$ of arbitrary genus with common dual $\Sigma^*$. By  \cite[Lemma $9.3$]{dave:lightbulb}, $\Sigma^*$ can be used to obtain disjointly embedded disks away from $\Sigma \cup \Sigma^*$ spanning a generating set for $H_1(\Sigma)$. Taking a regular neighborhood of these disks together with $\Sigma \cup \Sigma^*$, and then surgering $\Sigma^*$, yields an embedded $4$-ball $B \subset X^*$. Likewise, there is a $4$-ball $B' \subset X^*$ corresponding to $\Sigma'$. An argument analogous to the one above, using the $4$-balls $B, B' \subset X^*$, then produces a pairwise diffeomorphism $(X, \Sigma) \to (X, \Sigma')$. 
\end{remark}

\begin{remark}\label{iso}
Lemma \ref{diff} shows that the spheres $\Sigma_0$ and $\Sigma_n$ are equivalent (defined in \S0) in $X$ for all $n$: since the $\Sigma_n$ are homologous, the self-diffeomorphism given by the lemma induces the identity on $H_2(X)$. 
\end{remark}

\begin{theorem}\label{bounded}
The $2$-spheres $\Sigma_i$ and $\Sigma_j$ smoothly embedded in $X$ are topologically concordant in $X \times I$ if and only if $i \equiv j \mod 4$.   
\end{theorem}

\fig{210}{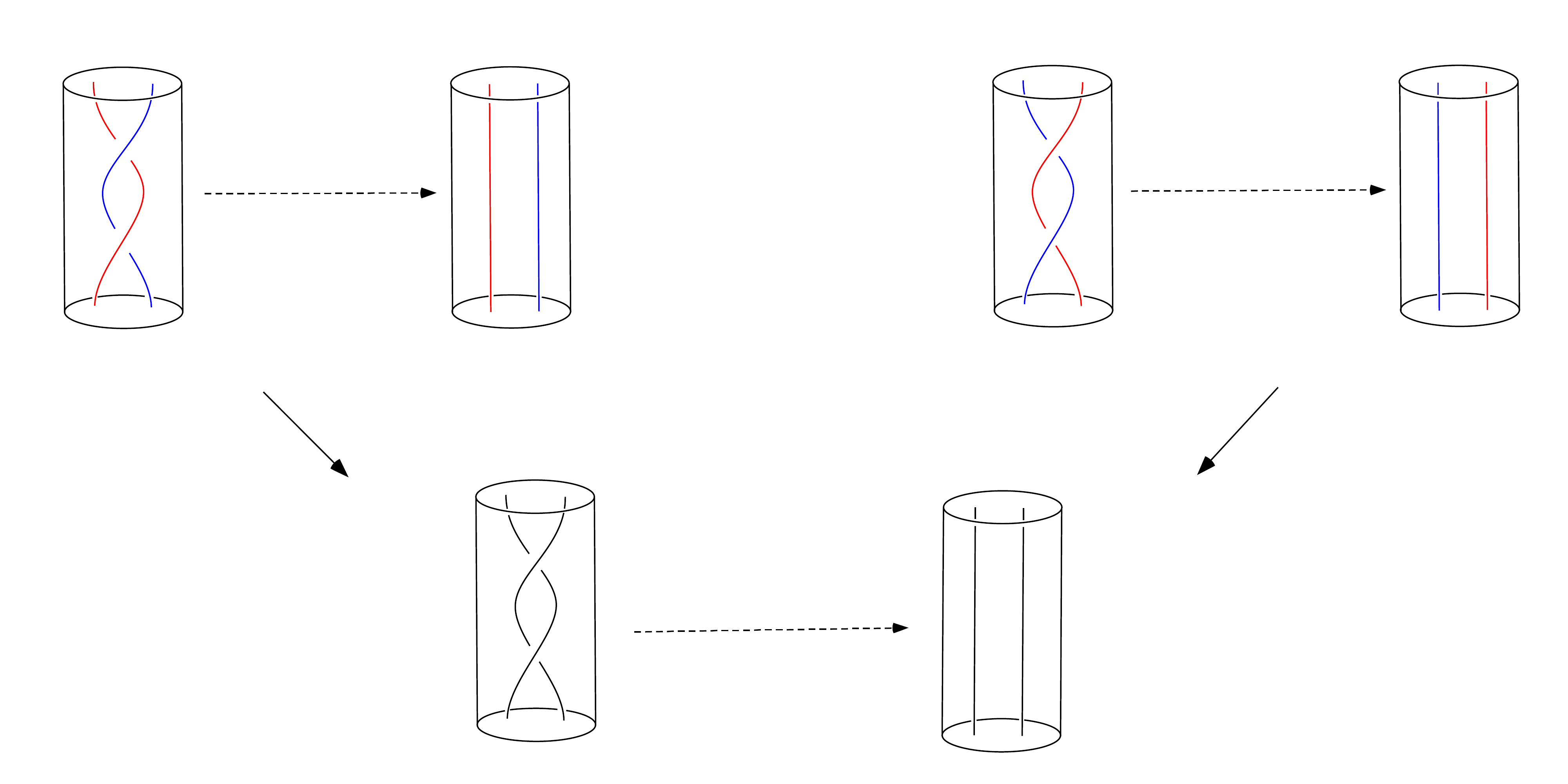}{
\put(-242,55){\small{Concordance $h$}}
\put(-228,46){\small{in $B \times I$}}
\put(-367,174){\small{Concordances}}
\put(-360,163){\small{in $B_0 \times I$}}
\put(-114,174){\small{Concordances}}
\put(-104,163){\small{in $B_\pi \times I$}}
\put(-81,88){$\pi$}
\put(-351,88){$\pi$}
\caption{The intersection of the $3$-ball $B=\{0\} \times D^2 \times I \subset X^{(1)}$ with $\Sigma_2$ (left bottom) and $\Sigma_0$ (right bottom), along with the intersection of its lifts $B_0$ and $B_\pi$ with both $\Sigma_{n,r}$ and $\Sigma_{n,b}$ for $n=0,2$.}
\label{lifts}}

\pf  By Theorem \ref{spherediff}, the spheres $\Sigma_i$ and $\Sigma_j$ are smoothly isotopic in $X$, and therefore smoothly concordant in $X \times I$, when $i \equiv j \mod 4$. On the other hand, as noted in the proof of Theorem \ref{spherediff}, the lifts of these 2-spheres to the universal $2$-fold cover $\widetilde X$ are homologous, and hence homotopic, if and only if $i$ and $j$ have the same parity. So, $\Sigma_i$ and $\Sigma_j$ are not even homotopic in $X$ if $i \not \equiv j \mod 2$. 

Therefore, it only remains to show that the homotopic spheres $\Sigma_0$ and $\Sigma_2$ are not topologically concordant in $X \times I$ (or equivalently that $\Sigma_1$ and $\Sigma_3$ are not concordant -- these two pairs of spheres are related by an automorphism of $X$). So, suppose to the contrary that there exists a topologically locally-flat concordance $h:S^2 \times I \to X \times I$ from $\Sigma_2$ to $\Sigma_0$. For each $\theta \in S^1$, the 2-sphere $\Sigma_n$ intersects the $3$-ball $\{\theta\} \times D^2 \times I \subset X^{(1)}$ in a properly embedded tangle that consists of a pair of arcs with endpoints where the attaching curves $\gamma_\pm$ of the 2-handles intersect the boundary of the $3$-ball, shown for $\Sigma_2$ and $\Sigma_0$ on the bottom of Figure \ref{lifts}. In general, the tangle corresponding to $\Sigma_n$ has $n$ half twists.

The universal (double) cover $\pi: \widetilde X \to X$ has deck transformation the involution $\tau$ induced by rotating the handlebody picture of $\widetilde X$ from Figure \ref{mfld} by angle $\pi$ about the dotted circle, as indicated in Figure \ref{lifthandles}. Let $\widetilde X^{(1)} \subset \widetilde X$ denote the union of the $0$-handle and the $1$-handle of $\widetilde X$, as with $X^{(1)} \subset X$ this is a copy of $S^1 \times D^2 \times [-1,1]$. The map $\pi$ winds $\widetilde X^{(1)}$ twice around $X^{(1)}$; intuitively, the lift ``unwraps" the generator of $\pi_1(X)$. So, for each $\theta \in S^1$, the $3$-balls $B_\theta =\{\theta\} \times D^2 \times I$ and $B_{\theta + \pi} =\{\theta + \pi\} \times D^2 \times I$ in $\widetilde X^{(1)}$ are both mapped to $\{2 \theta\} \times D^2 \times I \subset X^{(1)}$. Likewise, each 2-handle of $X$ lifts to a pair of 2-handles in $\widetilde X$. In particular, the 2-handles with framings $p$ and $q$ each lift to a pair of 2-handles, colored \emph{red} and \emph{blue}, that are interchanged by $\tau$ (see Figure \ref{lifthandles}). The attaching regions $S^1 \times D^2$ of these 2-handles intersect each 2-sphere $S_\theta= \partial B_\theta$ in red and blue co-cores $D^\pm_{r, \theta}$ and $D^\pm_{b,\theta}$ (one for each red/blue $2$-handle attached along a lift of $\gamma_\pm$).

Each $2$-sphere $\Sigma_n$ lifts to two spheres $\Sigma_{n,r}$ and $\Sigma_{n,b}$ in $\widetilde X$; again there is both a red and blue lift. Likewise, by the homotopy lifting property, the concordance $h$ lifts to a pair of concordances $h_r: S^2 \times I \to \widetilde X \times I$ from $\Sigma_{2,r}$ to $\Sigma_{0,r}$ and $h_b: S^2 \times I \to \widetilde X \times I$ from $\Sigma_{2,b}$ to $\Sigma_{0,b}$. Since the cover $\widetilde X$ is simply-connected, and both the red spheres and the blue spheres have a common dual (the two lifts of the common dual of the $\Sigma_n$), it follows from \cite[Theorem 10.1]{dave:lightbulb} that there is a smooth ambient isotopy from $\Sigma_{2,r} \sqcup \Sigma_{2,b}$ to $\Sigma_{0,r} \sqcup \Sigma_{0,b}$ in $\widetilde X$. However, we shall show that this is not ``equivariantly" true, i.e. no such isotopy can be the lift of an isotopy from $\Sigma_2$ to $\Sigma_0$ in $X$.

\fig{280}{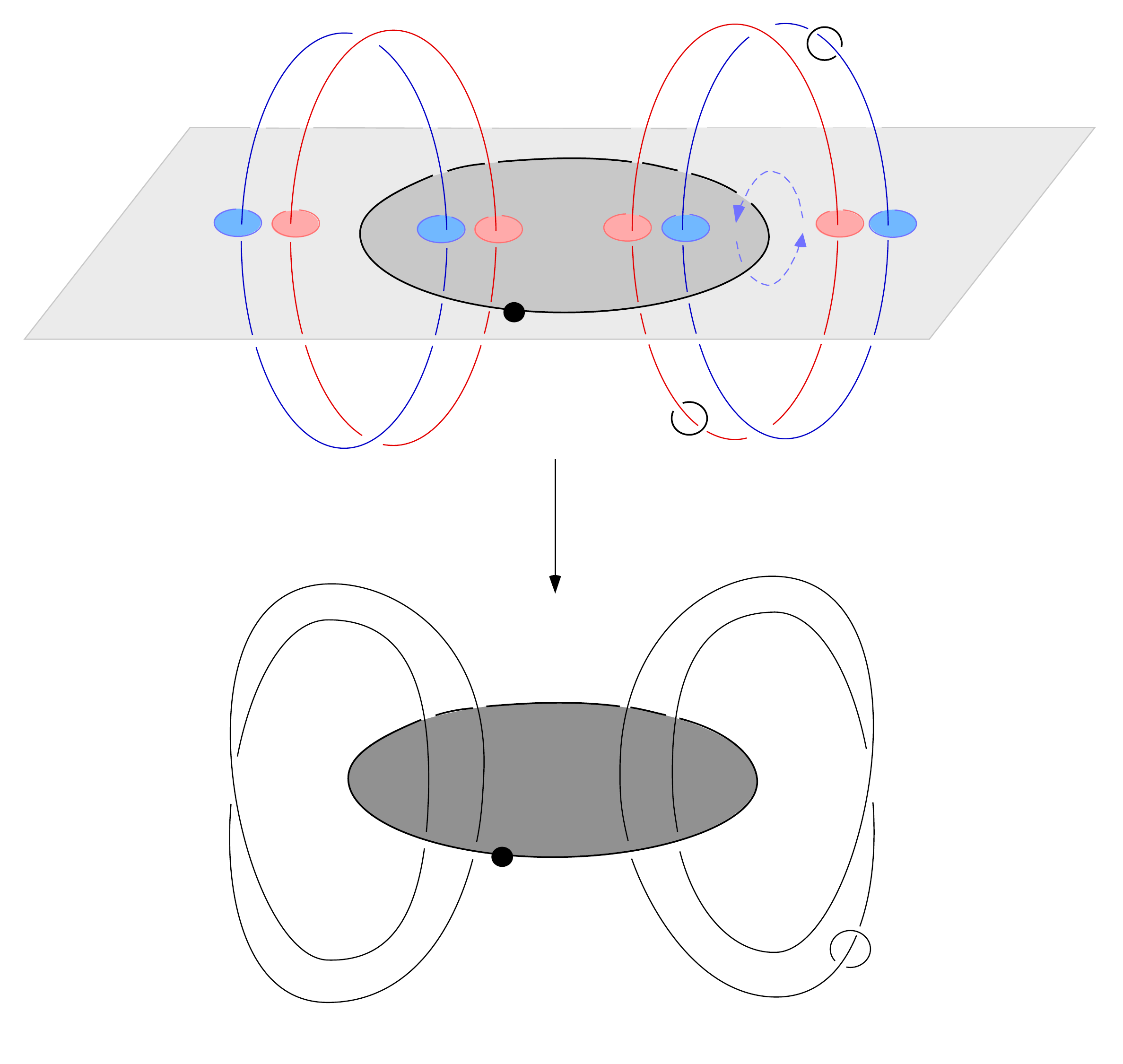}{
\put(-158,70){$\partial B$}
\put(-153,214){$S_0$}
\put(-260,200){$S_\pi$}
\put(-160,140){$\pi$}
\put(-246,100){$p$}
\put(-63,24){$0$}
\put(-69,266){$0$}
\put(-128,164){$0$}
\put(-70,164){$q-1$}
\put(-145,264){$q-1$}
\put(-174,164){$p+1$}
\put(-248,264){$p+1$}
\put(-60,100){$q$}
\put(-90,236){$\tau$}
\caption{The double cover $\pi: \widetilde X \to X$, with deck transformation $\tau$.}
\label{lifthandles}}

Let $\mathcal{T}\subset \widetilde X \times I$ denote the disjoint union $h_r(S^2 \times I) \sqcup h_b(S^2 \times I)$. Indeed, we may extend $X$ by a collar on its boundary in order for $\mathcal T \subset \interior(\widetilde X) \times I$. The concordance $h$ can then be isotoped rel boundary so that $\mathcal{T}$ is topologically transverse\foot{See Quinn \protect\cite[Theorem 2.4.1]{quinn}, and Freedman and Quinn \protect\cite[Theorem 9.5A]{freedman-quinn:4-manifolds}.} to the following subsets of $\widetilde X \times I$:

\begin{enumerate}
\item the cylinders $D^\pm_{r,\theta} \times I$ and $D^\pm_{b,\theta} \times I$ for $\theta \in \{0, \pi\}$, 
\item the $4$-balls $\Delta_r \times I$ and $\Delta_b \times I$, where $\Delta_r=\bigcup_{\theta \in [0,\pi]} D^+_{r,\theta}$ and $\Delta_b=\bigcup_{\theta \in [0,\pi]} D^+_{b,\theta}$ are each an interval's worth of co-cores in $\widetilde X$, and
\item the $4$-balls $B_0 \times I$ and  $B_\pi \times I$. 
\end{enumerate}  

To acheive this, we first arrange for $\mathcal T$ to be transverse to a collar neighborhood of the cylinders in $(1)$; that we can do so is immediate from \cite[Theorem 9.5A]{freedman-quinn:4-manifolds} since each cylinder is properly embedded in $\widetilde X \times I$. Note that at this point, $\mathcal T$ is transverse to the boundaries of the submanifolds in $(2)$ and $(3)$, which intersect $\interior(\widetilde X) \times I$ (and hence $\mathcal T$) exactly in the interior of the cylinders from $(1)$. It follows that $\mathcal T$ can also be made transverse to the submanifolds in $(2)$ and $(3)$, leaving the transverse intersection with the cylinders fixed.

\fig{200}{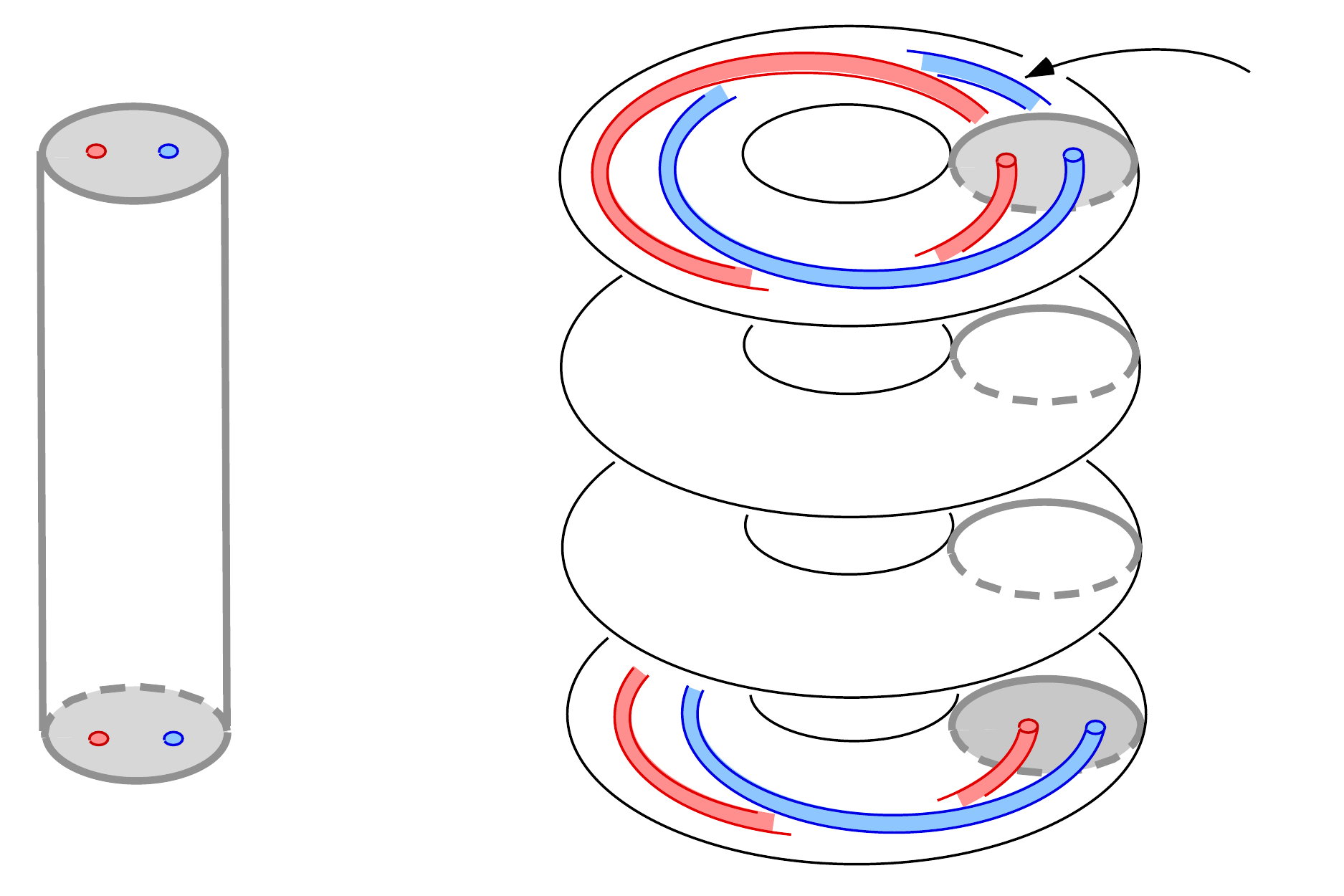}{
\put(-40,170){$2$-handle attaching regions}
\put(-215,94){$\subset$}
\put(-274,10){$S_0$}
\put(-155,-10){$\widetilde X^{(1)}=S^1 \times D^2 \times I$}
\caption{The sphere $S_0 \subset \widetilde X$. Note that the only portions of $S_0$ contained in the \emph{interior} of $\widetilde X$ are the disks $D^\pm_{r,0}$ and $D^\pm_{b,0}$ (the blue and red disks shown on $S_0$) where $S_0$ intersects the attaching regions of the $2$-handles.}
\label{lastpic}}

After these adjustments, $\mathcal{T}$ intersects the $4$-ball $B_0 \times I \subset \widetilde  X \times I$ in a properly embedded surface with both red and blue components, whose boundary is a link $\mathcal{L}$ embedded in $\partial(B_0 \times I) =S^3$. Note that a choice of orientation of the sphere $\Sigma_0$ induces an orientation of the concordance, and so also of $\mathcal{L}$. The linking number $\lk(\mathcal{L})$ between the red and blue components of $\mathcal L$ is equal to zero, since these components of the link bound disjointly embedded surfaces in the $4$-ball. The remainder of the proof will be spent contradicting this, as a further analysis of the link $\mathcal L$ shows that $\lk(\mathcal{L})$ must in fact also be odd. To begin, decompose the link $\mathcal{L}$ into the three oriented tangles illustrated in Figure \ref{tangles}: 
\begin{enumerate}[(i)]
\item a tangle $(\Sigma_{2,r} \cup \Sigma_{2,b}) \cap (B_0 \times \{0\})$ in the ``inner" 3-ball, 
\item the \emph{reverse mirror image} of the tangle $(\Sigma_{0,r} \cup \Sigma_{0,b}) \cap (B_0 \times \{1\})$ in the ``outer" 3-ball, 
\item and a tangle $\mathcal{T} \cap (S_0 \times I)$ connecting the endpoints of the tangles in (i) and (ii). 
\end{enumerate}

\fig{226}{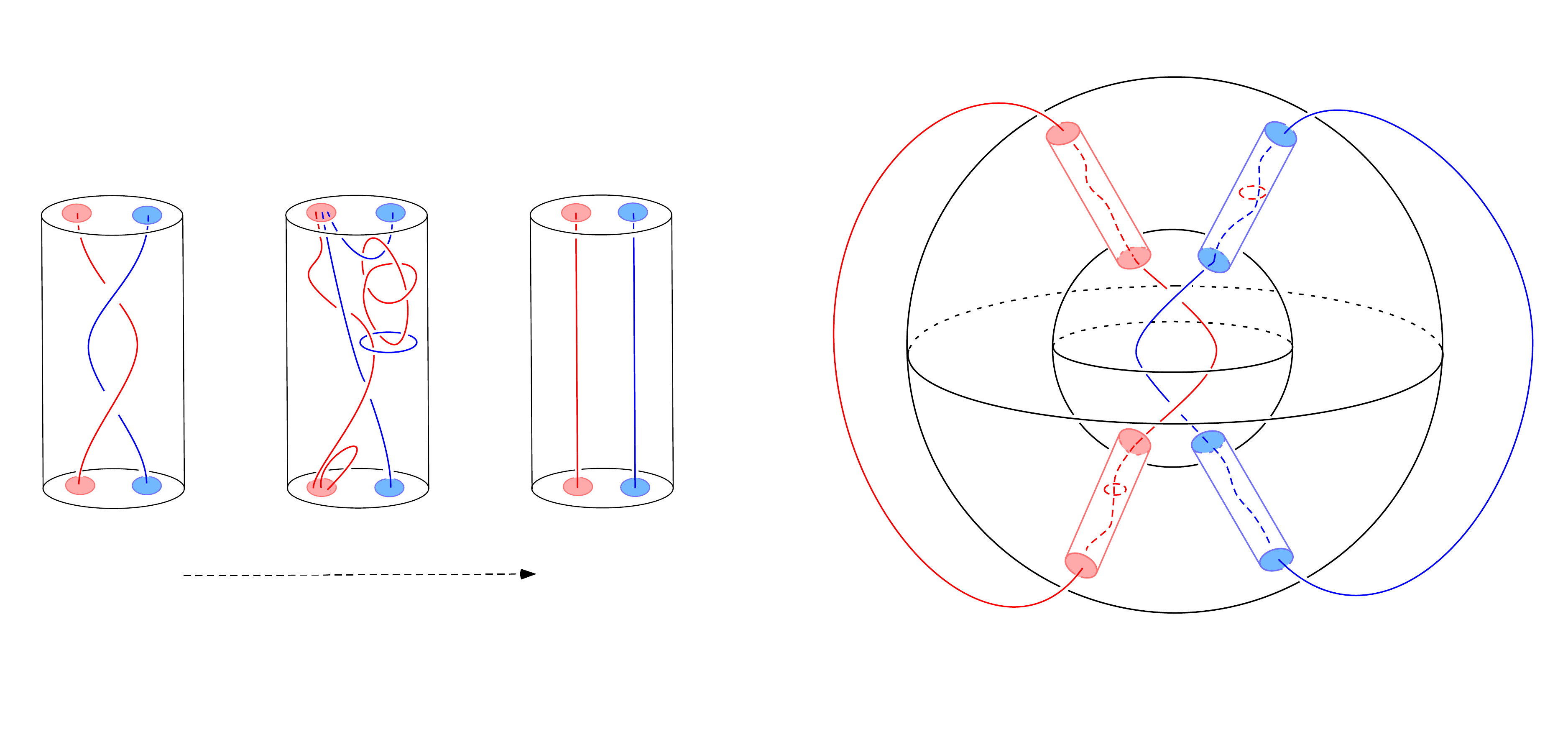}{
\put(-399,26){\small{$\mathcal{T} \cap (B_0 \times I)$}}
\put(-255,103){$=$}
\put(-254,111){\small{$\partial$}}
\put(-460,45){$t=0$}
\put(-307,45){$t=1$}
\put(-89,156){$T^+_{b,0}$}
\put(-169,73){$T^-_{r,0}$}
\put(-169,156){$T^+_{r,0}$}
\put(-89,73){$T^-_{b,0}$}
\put(-154,21){$\mathcal{L} \subset \partial(B_0 \times I)$}
\caption{The link $\mathcal{L}$, as the boundary of $\mathcal{T} \cap (B_0 \times I)$ in $\partial(B_0 \times I)$.}
\label{tangles}}

In fact, $\mathcal{T} \cap (S_0 \times I)$ may be decomposed further into the disjoint union of the four tangles $T^{\pm}_{r,0}=\mathcal{T} \cap (D^\pm_{r,0} \times I)$ and $T^\pm_{b,0}=\mathcal{T} \cap (D^\pm_{b,0} \times I)$ where $\mathcal{T}$ intersects each of the cylinders $D^\pm_{r,0} \times I$ and $ D^\pm_{b,0} \times I$; recall that these cylinders are exactly the regions of $S_0 \times I$ contained in $\interior(\widetilde X) \times I$, and hence the only portions of $S_0 \times I$ that intersect $\mathcal T$ (see Figure \ref{lastpic} for more explanation). Each of these tangles consists of a single properly embedded red (or blue) arc, and possibly many closed red and blue components. Examples of the tangles $T^\pm_{r,0}$ and $T^\pm_{b,0}$ are given in Figure \ref{tangles}.

Define tangles $T^\pm_{r,\pi}$ and $T^\pm_{b,\pi}$ contained in the cylinders $D^\pm_{r,\pi} \times I$ and $D^\pm_{b,\pi} \times I$ similarly. For each angle $\theta \in \{0,\pi\}$, connecting the endpoints of the tangles $T^\pm_{r,\theta}$ and $T^\pm_{b,\theta}$ by embedded, oriented arcs pushed into the interior of the cylinder from its boundary gives unique oriented links in the $3$-ball $\widehat T^\pm_{r,\theta}$ and $\widehat T^\pm_{b,\theta}$ which we call the ``closures" of the tangles. Let $\lk(\widehat T^\pm_{r,\theta})$ and $\lk(\widehat T^\pm_{b,\theta})$ denote the linking numbers between the red and blue components of these links. It is clear from Figure \ref{tangles} that 
\begin{align} \label{math} 
\lk(\mathcal{L}) \ = \ \lk(\widehat T^+_{r,0}) + \lk(\widehat T^+_{b,0})+ \lk(\widehat T^-_{r,0}) + \lk(\widehat T^-_{b,0}) + 1
\end{align} 
 

\smallskip

\noindent {\bf Claim 1:} $\lk(\widehat T^\pm_{r,0})= \lk(\widehat T^\pm_{b,0})$ 

\pf The involution $\tau \times \id: \widetilde X \times I \to \widetilde X \times I$ maps the tangle $T^+_{r,0} \subset D^+_{r,0} \times I$ to the tangle $T^+_{b,\pi} \subset D^+_{b,\pi} \times I$. Since both $(\tau \times \id) \circ h_r=h_b$ and $(\tau \times \id) \circ h_b=h_r$, the red components of the tangle $T^+_r$ are mapped to the blue components of $T^+_{b,\pi}$, and vice versa. Thus $\lk(\widehat T^+_{r,0})=\lk(\widehat T^+_{b,\pi})$.

Since the disks $D^+_{b,0}$ and $D^+_{b,\pi}$ are co-cores of the same (blue) 2-handle of $\widetilde X$, they are connected by the interval's worth of co-cores $\Delta_b=\bigcup_{\theta \in [0,\pi]} D^+_{b,\theta}$ within its attaching region $S^1 \times D^2$. The cylinder $\Delta_b$ is shaded blue in Figure \ref{s0bdd}. Recall that we initially arranged for $\mathcal{T}$  to intersect $\Delta_b \times I \subset \widetilde X \times I$ transversally. Therefore, $\mathcal{T} \pitchfork (\Delta_b \times I)$ is an embedded cobordism between the tangles $T^+_{b,0}$ and $T^+_{b,\pi}$ in $B^3 \times I$. Indeed, since both $\Sigma_{0,b}$ and $\Sigma_{2,b}$ intersect each of the co-cores $D^+_{b,\theta}$ in a single point for all $\theta$, this cobordism is a product on the boundary. 

Therefore, it is equivalent to say that the links  $\widehat T^+_{b,0}$ and $\widehat T^+_{b,\pi}$ are cobordant in $B^3 \times I$, and so $\lk(\widehat T^+_{b,0}) =\lk(\widehat T^+_{b,\pi})$. It follows that $\lk(\widehat T^+_{b,0})=\lk(\widehat T^+_{r,0})$, as desired. The argument that $\lk(\widehat T^-_{r,0})= \lk(\widehat T^-_{b,0})$ is analogous, completing the proof of the claim. 

\fig{280}{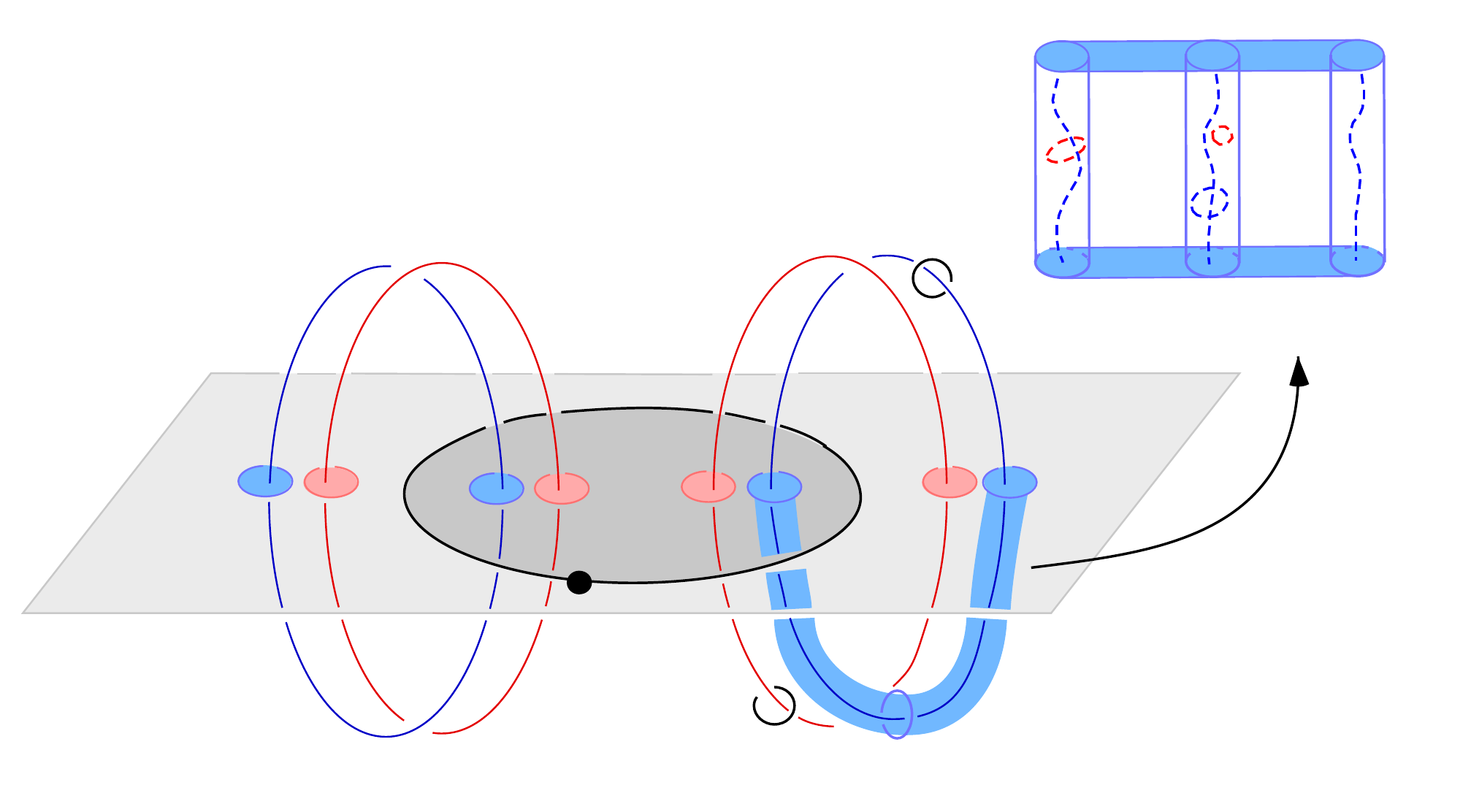}{
\put(-106,274){$\Delta_b \times I$}
\put(-68,104){$\times I$}
\put(-148,108){\small{$D^+_{b,\pi}$}}
\put(-231,108){\small{$D^+_{b,0}$}}
\put(-206,15){\small{$D^+_{b,\theta}$}}
\put(-207,108){\small{$D^+_{r,\pi}$}}
\put(-291,106){\small{$D^+_{r,0}$}}
\put(-148,172){\small{$T^+_{b,0}$}}
\put(-46,172){\small{$T^+_{b,\pi}$}}
\caption{The cobordism in $\Delta_b \times I$ from $T^+_{b,0}$ to $T^+_{b,\pi}$.}
\label{s0bdd}}

\smallskip

The claim, in conjunction with Equation (\ref{math}), implies that $ \lk(\mathcal{L}) \equiv \ 1 \pmod2$. However, this contradicts the fact that $\lk(\mathcal{L})=0$. Therefore, there can be no such concordance $h$. \qed 
\smallbreak 

This argument can be extended to produce a similar family of homotopic but non-isotopic smoothly embedded spheres in a \emph{closed} $4$-manifold; namely, the $2$-spheres $\Sigma_n \subset X$ sitting naturally in the closed double $\mathcal{D}X$, illustrated in Figure \ref{double}.  The double $\mathcal{D}X$ can be built from $X$ by adding a $(4-k)$-handle, thought of as an ``upside down" $k$-handle, for each $k$-handle of $X$. In particular, first attach a 2-handle along the boundary co-core of each 2-handle of $X$ (diagramatically, this amounts to attaching a $2$-handle along the meridian of the attaching circle of the 2-handle, using the zero framing). The resulting manifold is $\mathcal{D}X- (S^1 \times B^3)$, with boundary $S^1 \times S^2$. It follows that attaching one 3-handle along $\{\pt\} \times S^2$, and then one 4-handle, yeilds the double. Figure \ref{double2} shows how to picture the attaching sphere $\{\pt\} \times S^2$ when $S^1 \times S^2$ is drawn as the boundary of $\mathcal{D}X- (S^1 \times B^3)$. 

Since all of the extra 2-handles added to $X$ in order to form the double are attached along curves that are null-homotopic in $X$, there are isomorphisms $\pi_1(\mathcal{D}X) \cong \pi_1(X) \cong \mathbb{Z}_2$. The handlebody structure for $\mathcal{D}X$ induces a relative one for $\mathcal{D}X \times I$ (i.e. built from a collar of its boundary) whose $k$-handles correspond to the $(k-1)$-handles of $\mathcal{D}X$. 

\fig{320}{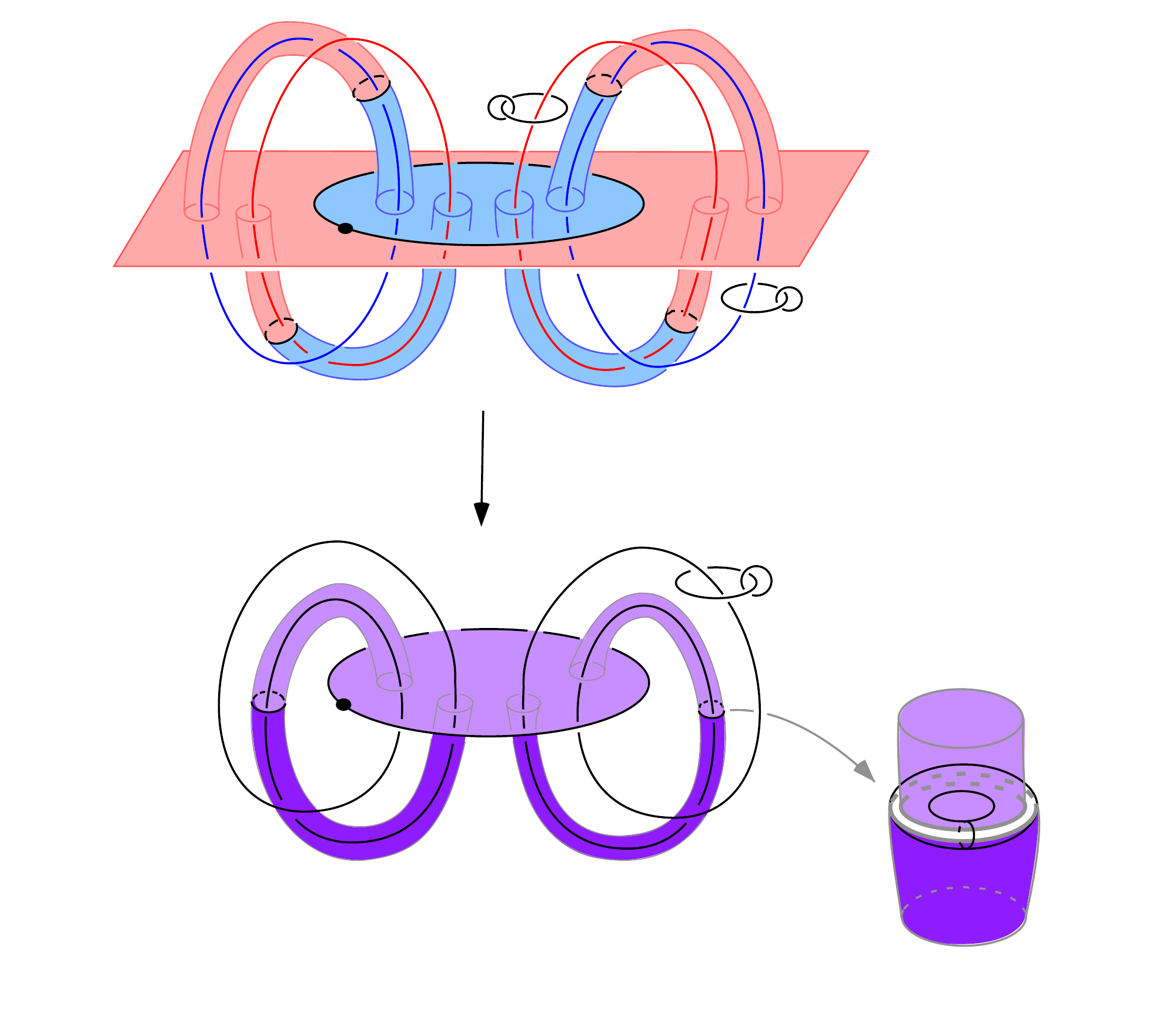}{
\put(-210,175){$\pi$}
\put(-274,96){\small $0$}
\put(-156,94){\small $0$}
\put(-142,146){\small $0$}
\put(-123,136){\small $0$}
\put(-306,112){\small $p$}
\put(-127,112){\small $q$}
\put(-150,200){\small $q-1$}
\put(-214,300){\small $q-1$}
\put(-236,200){\small $p+1$}
\put(-323,300){\small $p+1$}
\put(-30,40){\small{Surface passes twice}}
\put(-30,30){\small{over the 2-handle}}
\put(-30,20){\small{(with same sign).}}
\caption{The double $\mathcal{D}X$ (bottom) and its universal (2-fold) cover $\widetilde{\mathcal{D}X}$ (top). The attaching spheres for the 3-handles are pictured as punctured surfaces going over the 2-handles: the red and blue attaching 2-spheres are the lifts of the purple. The two sides of the purple surface are distinguished by their shading. Note that all unlabelled curves in the diagram of the cover are zero framed.}
\label{double}}

\fig{175}{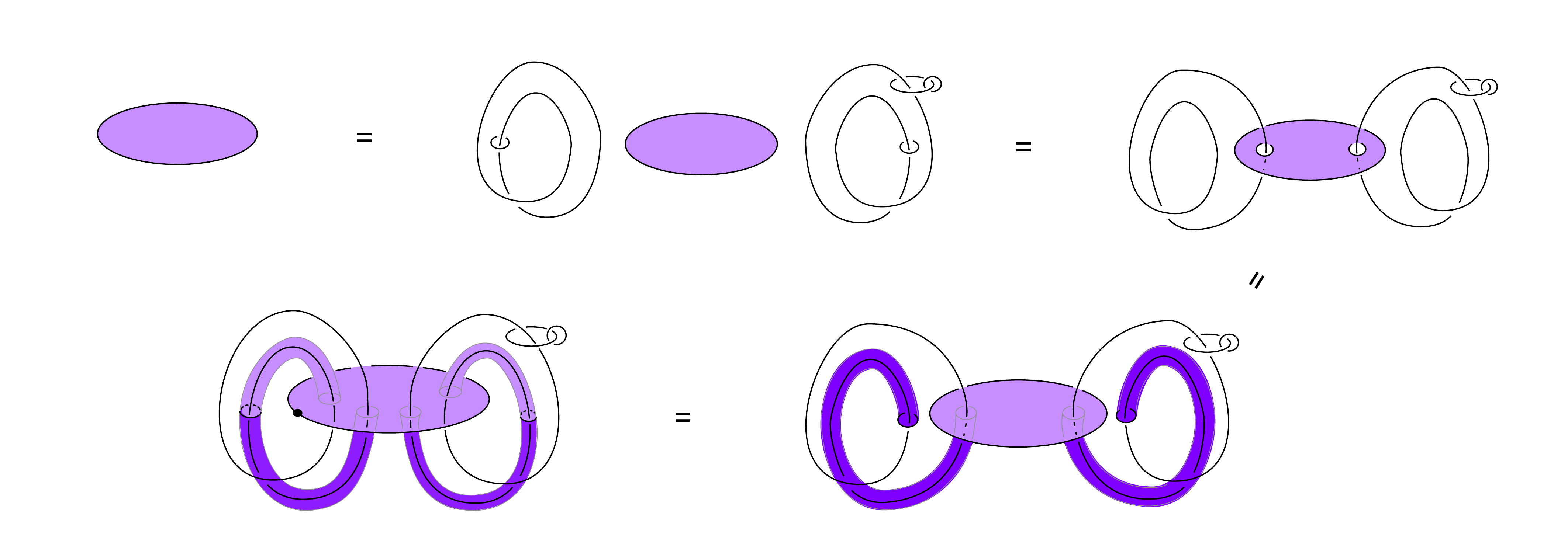}{
\put(-104,54){\small $q$}
\put(-434,54){\small $p$}
\put(-247,54){\small $p$}
\put(-316,54){\small $q$}
\put(-352,136){\small $p$}
\put(-200,136){\small $q$}
\put(-144,136){\small $p$}
\put(-24,136){\small $q$}
\put(-282,48){\small $\partial$}
\caption{The attaching 2-sphere (shaded purple)  in $S^1 \times S^2$ for the 3-handle of $\mathcal DX$. The first four figures in the sequence are $3$-manifolds diffeomorphic to $S^1 \times S^2$. All surgery coefficients are zero, unless specified otherwise. The last figure (bottom left) is the 4-manifold $\mathcal DX-(S^1 \times B^3)$, with the sphere in its boundary.}
\label{double2}}

\begin{theorem}\label{closed}
The $2$-spheres $\Sigma_i$ and $\Sigma_j$ smoothly embedded in $\mathcal{D}X$ are topologically concordant in $\mathcal{D}X \times I$ if and only if $i \equiv j \mod 4$.  
\end{theorem}

\pf As in the bounded case, by Theorem \ref{spherediff} the spheres $\Sigma_i$ and $\Sigma_j$ are smoothly isotopic in $X \subset \mathcal{D}X$, and therefore smoothly concordant in $\mathcal{D}X \times I$, when $i \equiv j \mod 4$, and are homotopic if and only if $i \equiv j \mod 2$. It again remains to show that $\Sigma_0$ and $\Sigma_2$ are not concordant. Suppose to the contrary that there exists a topological concordance $h:S^2 \times I \to \mathcal{D}X \times I$ from $\Sigma_2$ to $\Sigma_0$. Then $h$ lifts to a pair of red and blue concordances $h_r: S^2 \times I \to \widetilde{\mathcal{D}X} \times I$ and $h_b: S^2 \times I \to \widetilde{\mathcal{D}X} \times I$  whose images are disjointly embedded. 

We borrow much of the previous notation. Since $\widetilde X \subset \widetilde{\mathcal DX}$, the $3$-balls $B_\theta =\{\theta\} \times D^2 \times I$ and their boundary 2-spheres $S_\theta=\partial B_{\theta}$ are also embedded in $\widetilde{\mathcal DX}$. Following the previous proof, let the disks $D^\pm_{r,\theta}, D^\pm_{b,\theta} \subset S_\theta$ denote the co-cores where the attaching regions of the red and blue 2-handles of $\widetilde{\mathcal DX}$ intersect $S_\theta$. In this case, there is also a pair of annuli $A^\pm \subset S_0$ where the attaching regions of the 3-handles of $\widetilde{\mathcal DX}$ (after a perturbation) intersect $S_0$. It can be seen from Figure \ref{3handle} that $S_0$ intersects only the attaching region of the blue $3$-handle. Therefore, since they are disjoint, the annuli $A^+$ and $A^-$ are connected by an interval's worth of annuli $\Delta= (S^1 \times I) \times I$ in the attaching region of the blue $3$-handle.

Let $\mathcal{T}\subset \widetilde{\mathcal DX} \times I$ denote the disjoint union $h_r(S^2 \times I) \sqcup h_b(S^2 \times I)$. By general position, the image of the concordance $h$ can initially be arranged to avoid the $5$-handle of $\mathcal{D}X \times I$ (i.e. the product of the $4$-handle of $\mathcal DX$ with the interval $I$), so that $\mathcal {T}$ is embedded away from the $5$-handle of $\widetilde{\mathcal DX} \times I$. Therefore, the concordance $h$ may once again be isotoped rel boundary to arrange for $\mathcal{T}$ to be topologically transverse to the submanifolds $(1)$-$(3)$ listed in the previous proof, now thought of as subsets of $\widetilde{\mathcal DX} \times I$. In addition, we can arrange for $\mathcal{T}$ to be topologically transverse to the properly embedded solid tori $A^\pm \times I \subset \widetilde{\mathcal DX} \times I$, as well as $\Delta \times I \subset \widetilde{\mathcal DX} \times I$.

The first of these tranversality assumptions imply that $\mathcal{T} \pitchfork (B_0 \times I)$ is a properly embedded surface in the $4$-ball $B_0 \times I$ with boundary a link $\mathcal{L}$ with red and blue components in $\partial(B_0 \times I)=S^3$. Therefore, the linking number $\lk(\mathcal L)$ between the red and blue components is equal to zero. This fact, as in the previous proof, will be the source of our contradiction. 

Once more, decompose $\mathcal{L}$ into tangles in the 3-balls $B_0 \times \{0\}$ and $-B_0 \times \{1\}$ that are connected by a tangle in $S_0 \times I$. The tangle in $S_0 \times I$ is now composed not only of the bicolored tangles $T^\pm_{r,0} \subset D^\pm_{r,0} \times I$ and $T^\pm_{b,0} \subset D^\pm_{b,0} \times I$ (defined identically to those in the bounded case -- refer to Figure \ref{tangles}) but also two additional tangles $C^\pm = \mathcal{T} \cap (A^\pm \times I)$ shown in Figure \ref{s0closed}.  In fact the $C^\pm$ are links, since the red and blue 2-spheres at the ends of each concordance are disjoint from the attaching sphere of the blue 3-handle. 

\fig{135}{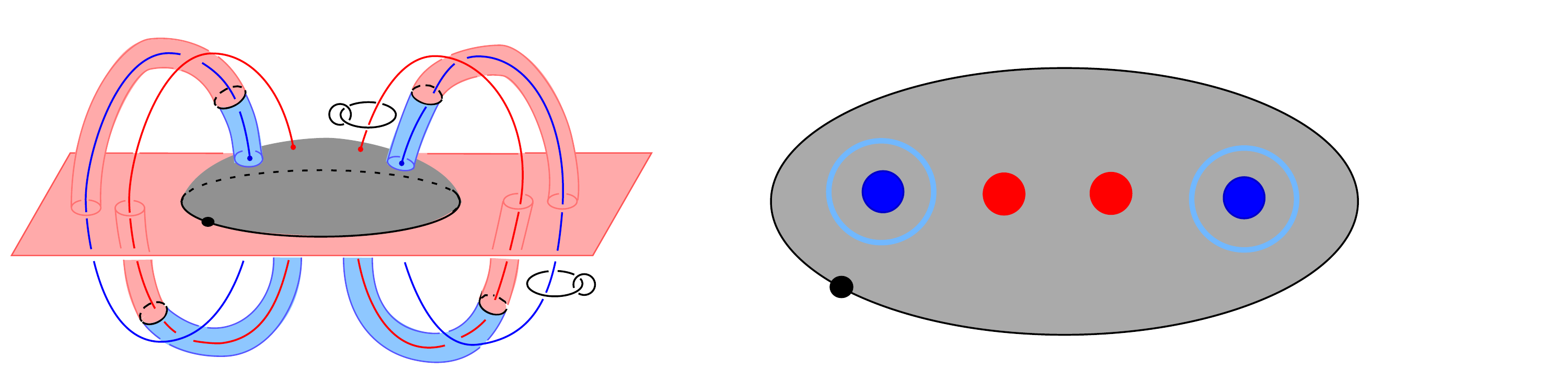}{
\put(-434,60){$S_0$}
\put(-177,28){$S_0$}
\put(-239,37){$A^+$}
\put(-120,36){$A^-$}
\caption{The intersection of the attaching regions of the $2$-handles and $3$-handles of $\widetilde{\mathcal{D}X}$ with the $2$-sphere $S_0$.}
\label{3handle}}

As in the previous proof, let $\widehat T^\pm_{r,0}$ and $\widehat T^\pm_{b,0}$ be the ``closures" of the tangles $T^\pm_{r,0} \subset D^\pm_{r,0} \times I$ and $T^\pm_{b,0}$. Let $\lk(\widehat T^\pm_{r,\theta})$ and $\lk(\widehat T^\pm_{b,\theta})$ denote the linking numbers between the red and blue components of these links. In addition, let $\lk(C^\pm)$ denote the linking number between the red and blue components of $C^\pm$, thought of in the  $3$-sphere $\partial(B_0 \times I)$. Recall that $C^\pm \subset \mathcal{L}$ inherits an orientation once an orientation of the $2$-sphere at one end of the concordance is chosen. Let $w_r(C^\pm)$ and $w_b(C^\pm)$ denote the integral winding numbers of the red and blue components of $C^\pm$ around the $S^1$ factor of the solid tori $A^\pm \times I$, once an orientation on the $S^1$ factor is fixed. Now, it is clear from Figure \ref{s0closed} that 
\begin{align} \label{math2} 
\lk(\mathcal{L}) \ = \ \lambda^+ + \lambda^- + w + \lk(C^+) + \lk(C^-) + 1 
\end{align} 
where $\lambda^\pm= \lk(T^\pm_{r,0}) + \lk(T^\pm_{b,0})$ and $w = w_r(C^+) + w_r(C^-) +  w_b(C^+) + w_b(C^-)$.

\smallskip

\noindent {\bf Claim 2:} {\bf (a)} $\lk(C^+)=\lk(C^-)$\\
\indent  \hskip 41pt {\bf (b)} $w_r(C^+)=w_r(C^-)$ and $w_b(C^+)=w_b(C^-)$ 

\smallbreak 

\fig{260}{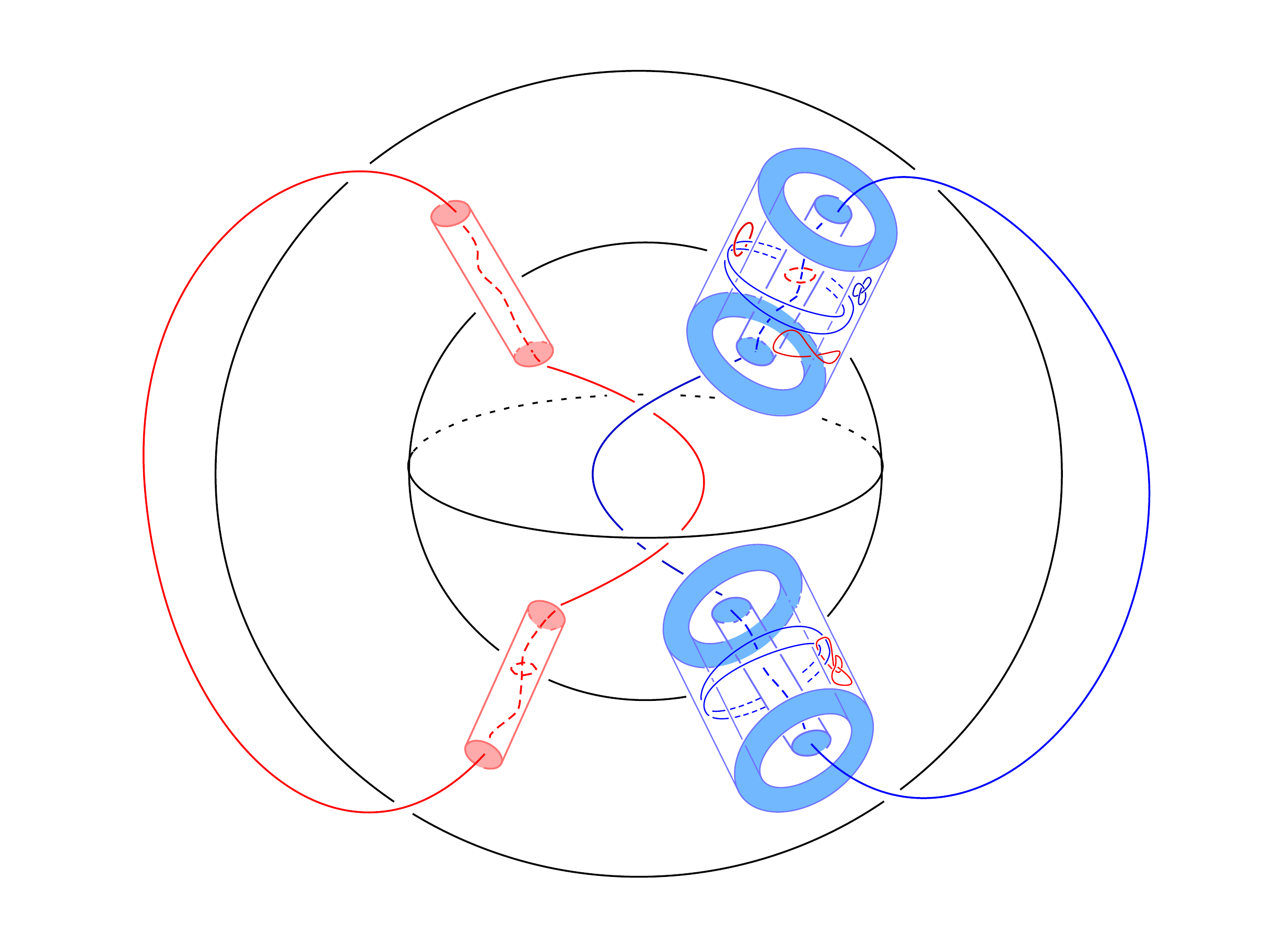}{
\put(-114,168){$C^+$}
\put(-117,70){$C^-$}
\caption{The link $\mathcal{L}$ in the closed case, embedded in the $3$-sphere $\partial(B_0 \times I)$ as the boundary of the surface $\mathcal{T} \cap (B_0 \times I)$.}
\label{s0closed}}

\noindent \emph{Proof}: The product $\Delta \times I$ can be thought of as a (relative) cobordism between the solid tori $A^+ \times I$ and $A^- \times I$. From this perspective, the transverse intersection $\mathcal{T} \pitchfork \Delta \times I$ is an embedded cobordism in $\Delta \times I$ from $C^+ \subset A^+ \times I$ to $C^- \subset A^- \times I$. This implies that the winding numbers of the red and blue sub-links of $C^\pm$ are preserved, as well as the linking number between the red and blue components, proving the claim. 
\smallbreak

Note that the argument from the bounded case still applies to show that $\lk(T^+_{r,0})= \lk(T^+_{b,0})$ and $\lk(T^-_{r,0})= \lk(T^-_{b,0})$.  Thus the sums $\lambda^\pm= \lk(T^\pm_{r,0}) + \lk(T^\pm_{b,0})$ are even, as is $w = w(C^+) + w(C^-)$. So Claims $(1)$ and $(2)$, in conjunction with Equation (\ref{math2}), imply that $\lk(\mathcal{L})  \equiv \ 1 \pmod2$. However, as this contradicts the fact that $\lk(\mathcal{L})=0$, there can be no such concordance $h$. \qed 

\smallskip

\begin{remark*} Theorem $1$ follows from Theorems \ref{bounded} and \ref{closed}, using the pair of $2$-spheres $\Sigma_0$ and $\Sigma_2$ that are equivalent by Remark \ref{iso}. To generate infinitely many $4$-manifolds as in Theorem 1, one can vary the framings $p$ and $q$ of the $2$-handles in either $X$ or its double. Alternatively, the arguments above apply to $4$-manifolds built from $X$ by adding arbitrary $2$-handles attached away from the $3$-ball $B_0 \subset X^{(1)}$. These include the connected sums of $X$ with multiple copies of $S^2 \times S^2$ or $\pm \cptwo$.  
\end{remark*}

We end with a question, as well as the proof of our second theorem. First, let $T_X$ denote the set of all (non-trivial) $2$-torsion elements of $\pi_1(X)$. Recall from Section $1$ that the crossed cycles of a smooth homotopy between embedded $2$-spheres  in a $4$-manifold $X$ correspond to a finite collection of (not necessarily distinct) elements in $T_X$. However, since Gabai's Generalized 4D Lightbulb Theorem shows that spheres related by a regular homotopy corresponding to a pair of identical elements in $T_X$ are smoothly isotopic, it is most relevant to record only the parity of each element's multiplicity; this gives an element of $\mathbb{Z}_2[T_X]$  (where $\mathbb{Z}_2[T_X]$ is considered equal to the trivial group when $T_X$ is empty). For example, the regular homotopy discussed in Appendix \ref{appendix} corresponds to the element $x \in \mathbb{Z}_2[T_X]$, where $x$ is the generator of $\pi_1(X) \cong \mathbb{Z}_2$.  

\smallskip

\noindent 
{\bf Question 1.}  Given a 2-sphere $S$ smoothly embedded in \emph{any} $4$-manifold $X$, and any non-zero element $h \in \mathbb{Z}_2[T_X]$, does there exist a $2$-sphere $T$ embedded in $X$ not isotopic to $S$ such that $S$ and $T$ are related by a regular homotopy corresponding to $h$? 

\begin{remark} In work since the release of the first version of this paper, Teichner and Schneiderman \cite{st} use the Freedman-Quinn invariant, defined in \cite[Theorem 10.5]{freedman-quinn:4-manifolds} and later corrected by Stong \cite{stong}, to show that Question 1 has an affirmative answer \emph{up to a distinguished subgroup $G \subset \mathbb{Z}_2[T_X]$} when $S$ has a dual sphere. In other words, given any $w \in \mathbb{Z}_2[T_X]$, there is a homotopy corresponding to $w$ taking $S$ to some sphere $T \subset X$. If $w \not \in G$, then one can use the Freedman-Quinn invariant to obstruct a topological isotopy between $S$ and $T$. However, if $w \in G$, then the homotopy can be replaced by one corresponding to $0 \in \mathbb{Z}_2[T_X]$. It then follows from Gabai's restated 4D Lightbulb Theorem \cite{dave:lightbulb} that $S$ and $T$ are smoothly isotopic. 
\end{remark} 

\smallskip

\noindent \emph{Proof of Theorem 2:} In \cite{morganszabo}, Morgan and Szab\'o exploit the chamber structure of the Seiberg-Witten invariant when $b_2^+=1$ to produce a family of rational surfaces $X$ (the proof here applies to any one of these) with h-cobordisms $W$ from $X$ to $X$ that are not diffeomorphic to $X \times I$. These cobordisms can be constructed explicitly as follows. Let $Z$ be the $5$-manifold obtained from $X\times I$ by attaching a $2$-handle along the boundary of an embedded disk $D$ in $X\times\{1\}$. Then $Z$ is a cobordism from $X$ to a $4$-manifold $X^\circ$ diffeomorphic to $X\cs\sss$, from now on considered literally equal. Indeed, attaching the $2$-handle using the appropriate framing, the sphere $A=S^2 \times \pt$ is the union of the disk $D$ capped off with a parallel copy of the core of the $2$-handle, and $A^*=\pt \times S^2$ is the belt sphere of the $2$-handle. The cobordism $W$ in the proof of \cite[Theorem $1.1$]{morganszabo} is formed by gluing $Z$ to $-Z$ along $X^\circ$ using a suitable (orientation preserving) automorphism $\phi: X^\circ \to X^\circ$. Using Wall \cite{wall:diffeomorphisms} to modify $\phi$ if necessary, and setting $B= \phi(A^*)$, we may assume that $A^* \cdot B  = \pm1$ and that $A$ and $B$ are homologous.  

Let $N \subset X^\circ$ be a regular neighborhood of the union $A \cup A^*$, diffeomorphic to a punctured copy of $S^2 \times S^2$. Surgering $A \subset N$ or $B \subset \phi(N)$ gives a corresponding 4-ball $D$ or $D'$ in $X$. By Palais \cite{palais}, there is an ambient isotopy of $X$ taking $D$ to $D'$. The end of this isotopy gives a diffeomorphism from $X-D$ to $X-D'$, or equivalently, from  $X^\circ -N$ to $ X^\circ-\phi(N)$. By Cerf \cite{cerf}, the restriction of this map to the boundary 3-sphere $\partial D$ may be isotoped to $\phi|_{\partial D}$ and then extended across $N$ to give an automorphism of $X^\circ$ carrying $A$ to $B$. Since the spheres $A$ and $B$ are homologous, this map acts like the identity on $H_2(X^\circ)$, and therefore the spheres $A$ and $B$ are equivalent. 

Since $A$ and $B$ are homologous and have simply-connected complement in $X^\circ$, they must also be topologically isotopic by Lee and Wilczy\'nski \cite{lw} (and more recently Sunukjian \cite{sunukjian:isotopy}). However, $A$ and $B$ are \emph{not} smoothly isotopic in $X^\circ$. If they were, then the 2-handle and the 3-handle of the cobordism $W$ would cancel geometrically; as $W$ is not smoothly a product, this is not the case. \qed

\begin{remark*}
It follows from Quinn \cite{quinn1} that the diffeomorphism constructed in the proof of Theorem $2$ above is \emph{topologically} isotopic to the identity, since it acts trivially on homology and the manifold $X^\circ$ is closed and simply-connected. The first examples of diffeomorphisms topologically but not smoothly isotopic to the identity were given by Ruberman \cite{danny1} using 1-parameter gauge theory; the proof above gives an alternate method of constructing such diffeomorphisms, in manifolds with $b_2^+=1$. Note that by virtue of using Morgan and Szab\'o's results from \cite{morganszabo}, our proof still relies on gauge theoretic techniques. 
\end{remark*}

\appendix 
\section{More on regular homotopy} \label{appendix} 
As promised, we proceed to describe an explicit regular homotopy from  $\Sigma_n$ to $\Sigma_{n+2}$ with a single crossed cycle corresponding to the generator of $\pi_1(X) \cong \mathbb{Z}_2$. As noted in Remark \ref{iso1}, the manifold $X$ is the union of the support of this regular homotopy with a neighborhood of the dual $2$-sphere, and so also provides a local model for such a homotopy in a larger manifold.
\smallskip

For each $n \in \mathbb{Z}$, after an isotopy of $\Sigma_n$, it may be assumed that $\Sigma_n$ intersects $(S^1 \times D^2) \times [0,\epsilon]$ in the product $\gamma \times [0,\epsilon]$, where $\gamma$ denotes the (2,1) torus curve pushed into the interior of the solid torus $S^1 \times D^2$, pictured on the right of Figure \ref{reghom4}. There is a regular homotopy $\rho_n: S^2 \times I \to X \times I$ from $\Sigma_n$ to $\Sigma_{n+2}$ supported in  $(S^1 \times D^2) \times [0,\epsilon]$ that consists of one finger move resulting in the immersed sphere $\Sigma_n'$ on the right of Figure \ref{reghom2}, followed by one Whitney move using the Whitney disk $W$ shown in blue. The embedded sphere $\Sigma_n'' \subset X$ obtained by performing this Whitney move is depicted in Figure \ref{reghom3}. In fact, there is a smooth isotopy supported in $(S^1 \times D^2) \times [0,\epsilon]$ taking $\Sigma_n''$ to the sphere $\Sigma_{n+2}$. 

\fig{340}{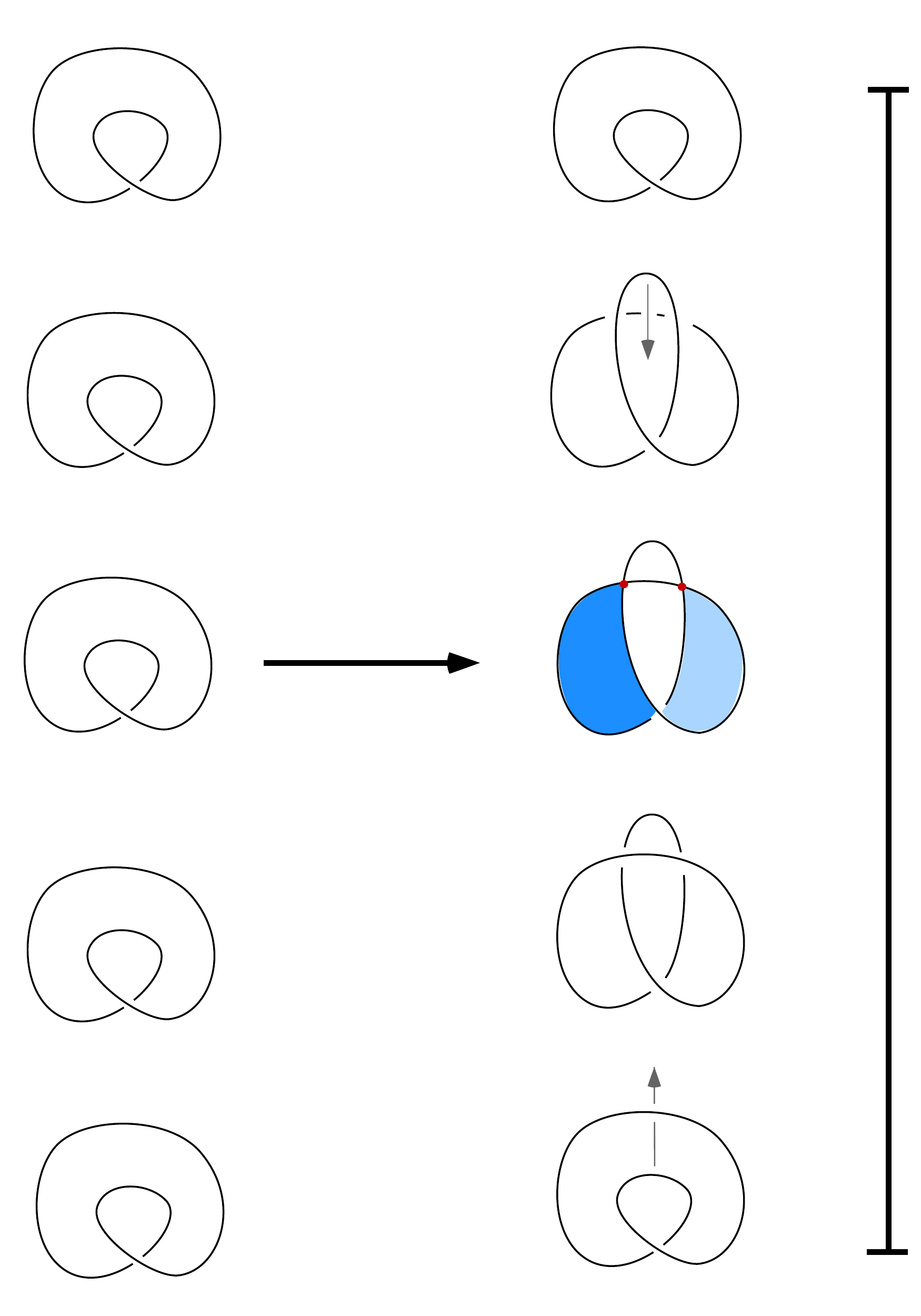}{
\put(-87,190){\small{$-$}}
\put(-60,190){\small{+}}
\put(-90,165){\small{$W$}}
\put(-170,150){\small{Finger move}}
\put(-155,140){\small{of $\rho_n$}}
\put(-223,-15){$S^1 \times D^2$}
\put(-88,-15){$S^1 \times D^2$}
\put(4,313){$\epsilon$}
\put(4,13){$0$}
\caption{The intersection of $\Sigma_n$ (left) and $\Sigma_n'$ (right) with $(S^1 \times D^2) \times [0,\epsilon]$. Each level set is a curve that should be thought of in the solid torus $(S^1 \times D^2) \times \{h\}$ for some height $h \in [0,\epsilon]$ (compare with Figure \ref{reghom1}). Also shown is the Whitney disk $W$ (in blue) used to guide the Whitney move in Figure \ref{reghom3}.}
\label{reghom2}}

\fig{250}{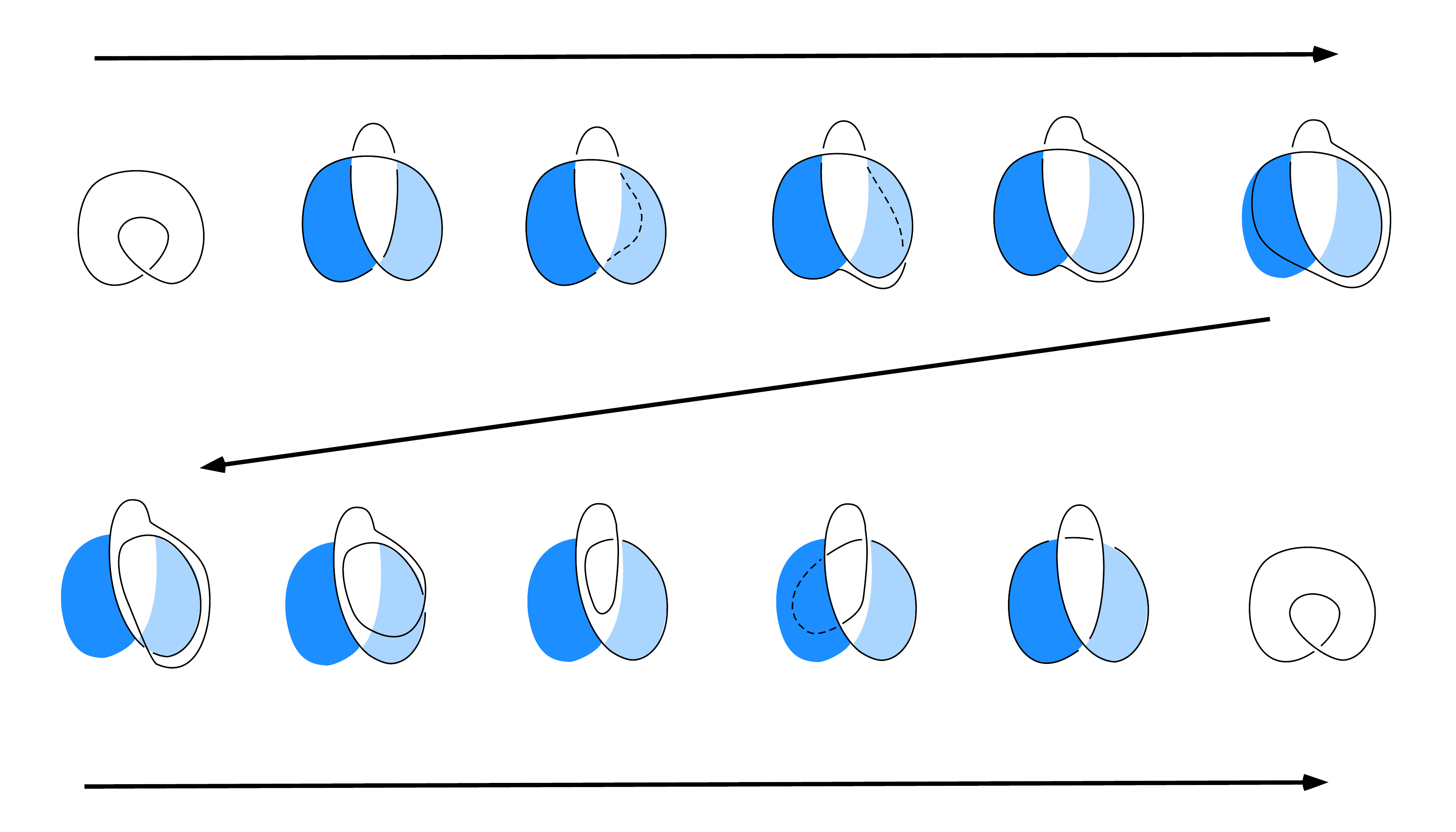}{
\put(-260,240){\small{Slide across front of $W$}}
\put(-260,20){\small{Slide across back of $W$}}
\caption{The intersection of $\Sigma_n''$ with $(S^1 \times D^2) \times [0,\epsilon]$. Each level set is a curve that should be thought of in the solid torus $(S^1 \times D^2) \times \{h\}$ for some height $h \in [0,\epsilon]$. As described in Figure \ref{reghom1}, one sheet of $\Sigma_n''$ slides first across the front and then the back of the Whitney disk $W$ (the two sides of the blue disk are distinguished by their shading).}
\label{reghom3}}

Figure \ref{reghom4} shows the single crossed cycle of $\rho_n$, which corresponds to the generator of $\pi_1(X) \cong \mathbb{Z}_2$. In fact, there is a regular homotopy in $X$ with $k$ crossed cycles from $\Sigma_n$ to $\Sigma_{n+2k}$ gotten by concatenating the $k$ regular homotopies $\rho_n, \rho_{n+2},\dots, \rho_{n+2k-2}$. It is then immediate from the restated version of the 4D Lightbulb Theorem that the 2-spheres $\Sigma_i$ and $\Sigma_j$ are smoothly ambiently isotopic in $X$ when $i \equiv j \mod 4$.

\fig{140}{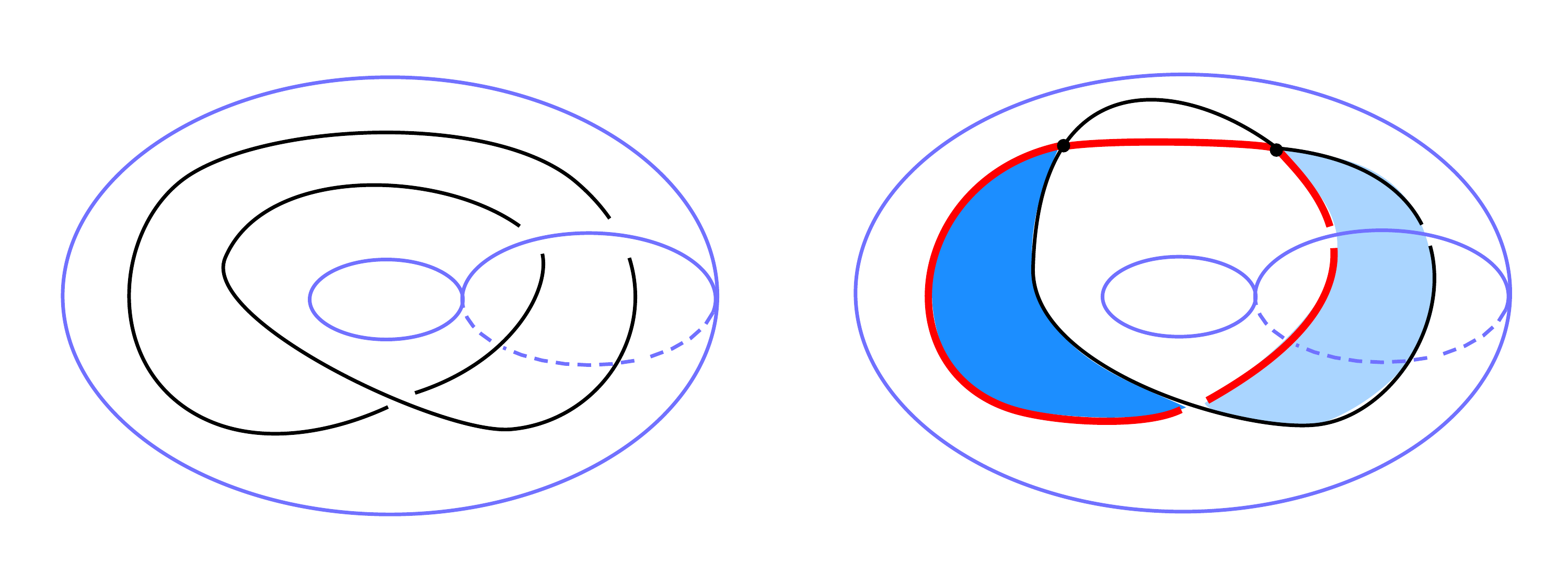}{
\put(-270,30){\small{$\gamma$}}
\put(-94,97){\small{$\eta$}}
\put(-144,66){\small{W}}
\caption{The curve $\gamma \subset S^1 \times D^2$ (left), and the curve $\eta$ of double points of the homotopy $\rho_n$ projected from $X \times I$ to $X$ (shown on the right, in red). Since $\eta$ is nontrivial in  $\pi_1(X) \cong \mathbb{Z}_2$ and has a single curve in its preimage, this is a crossed cycle corresponding to the generator of  $\pi_1(X) \cong \mathbb{Z}_2$.}
\label{reghom4}}

\end{document}